\documentclass[12pt,a4paper]{article}
\usepackage{amssymb}
\usepackage{amsmath,amsfonts,amsthm}

\usepackage{graphicx}
\usepackage{amsmath}

\usepackage{amsfonts}
\usepackage{amsthm,amscd}
\usepackage{graphicx}
\usepackage{amsmath}
\usepackage{amssymb}
\usepackage{amstext}

\newtheorem{theorem}{Theorem}
\newtheorem{corollary}[theorem]{Corollary}
\newtheorem{definition}[theorem]{Definition}
\newtheorem{example}[theorem]{Example}
\newtheorem{lemma}[theorem]{Lemma}

\newtheorem{proposition}[theorem]{Proposition}
\newtheorem{remark}[theorem]{Remark}

\title{Haar systems, KMS states on von Neumann algebras and $C^*$-algebras on   dynamically defined groupoids and Noncommutative Integration}

\author{G. G. de Castro,  A. O. Lopes  and G. Mantovani}

\begin{document}

\maketitle

\begin{abstract}

We analyse Haar systems associated to groupoids obtained by certain equivalence relations of   dynamical nature on sets like  $\{1,2,...,d\}^\mathbb{Z}$,  $\{1,2,...,d\}^\mathbb{N}$, $S^1\times S^1$, or $(S^1)^\mathbb{N}$, where $S^1$ is the unitary circle. We also describe properties of   transverse functions, quasi-invariant probabilities and KMS states for some examples  of  von Neumann algebras (and also $C^*$-Algebras) associated to these groupoids. We relate some of these KMS states with Gibbs states of Thermodynamic Formalism.
While presenting new results, we will also describe in detail several examples and basic results on the above topics. In other words it is also a survey paper.
Some known results on non-commutative integration are  presented, more precisely,  the relation of transverse measures, cocycles and quasi-invariant probabilities.

We describe the results in a language which is more familiar to people in  Dynamical
Systems.
Our intention is to study
Haar systems, quasi-invariant probabilities and von Neumann algebras as a topic on measure theory (intersected with ergodic theory) avoiding
questions of algebraic nature (which, of course, are also extremely  important).

\end{abstract}

\section{The groupoid associated to a partition} \label{gru}

\bigskip

We will analyze properties of Haar systems, quasi-invariant probabilities, transverse measures, $C^*$-algebras and KMS states related to Thermodynamic Formalism and Gibbs states. We will consider a specific particular setting where the groupoid will be defined by some natural  equivalence relations on the sets of the form $\{1,2,...,d\}^\mathbb{N}$ or $\{1,2,...,d\}^\mathbb{Z}$, $S^1\times S^1$, or $(S^1)^\mathbb{N}$. These  equivalence relations will be of dynamic origin.

We will denote by $X$ any one of the above sets.

The main point here is that we will use a notation which is more close to the one used on  Ergodic Theory and Thermodynamic Formalism.

 On section \ref{ktran} we introduce the concept of transverse functions associated to groupoids and Haar systems.

 On section \ref{quasi} we  consider modular functions and quasi-invariant probabilities on groupoids. In the end of this section we present a new result concerning a (non-)relation of the quasi-invariant probability  with the SBR probability of the generalized Baker map.

 On section \ref{von1} we consider a certain von Neumann algebra and the associated KMS states.  On proposition \ref{poa1} we present a new result concerning the relation between  probabilities satisfying the KMS property (quasi-invariant) and Gibbs (DLR) probabilities of Thermodynamic Formalism on the symbolic space $\{1,2,...,d\}^\mathbb{N}$ for a certain groupoid. Proposition \ref{KMS1} shows that the KMS probability is not unique on this case.

\cite{Hahn2}, \cite{FM1} and \cite{FM2} are the classical references on measured groupoids and von Neumann algebras. KMS states and $C^*$-algebras are described on \cite{Pana}.

 On  section  \ref{non} we present a natural expression - based on quasi-invariant probabilities - for the integration of a transverse function by a transverse measure. Some basic results on non-commutative integration (see \cite{Con} for a detailed description of the topic) are briefly described.

 On section \ref{C2} we present briefly the setting of $C^*$-algebras associated to groupoids on symbolic spaces. We present the well known and important concept of  approximately proper equivalence relation and its relation with the direct inductive limit topology (see \cite{Exel1}, \cite{Exel2}, \cite{Exel3} and \cite{Ren2}).

On section \ref{exaquasi} we present several examples of quasi-invariant probabilities for different kinds of groupoids and Haar systems.

On section \ref{ergo} we briefly describe ergodic concepts for groupoids and measures.

Results on $C^*$-algebras and KMS states from the point of view of  Thermodynamic Formalism   are presented in \cite{KumRen}, \cite{Ren2}, \cite{Put}, \cite{Thom1}, \cite{Thom2}, \cite{AHR}, \cite{EL1} and \cite{EL2}.

The paper \cite{Bis}  considers equivalence relations  and DLR probabilities for certain interactions on the symbolic space
$\{1,2,...d\}^\mathbb{Z}$ (not in $\{1,2,...d\}^\mathbb{N}$ like here).

Theorem 6.2.18 in Vol II of
\cite{Bra} and \cite{Araki} describe the relation between KMS states  and  Gibbs probabilities  for  interactions on certain spin lattices (on the one-dimensional case  corresponds to the space $\{1,2,..,d\}^\mathbb{Z}).$

We point out that Lecture 9 in \cite{GiaII} presents a brief introduction to $C^*$-Algebras and  non-commutative integration.

\medskip

We are indebted to Ruy Exel for many helpful discussions and useful comments during the procedure of writing  this paper. We also thanks Ali Tahzibi for some fruitful remarks.

\medskip

We denote   $\{1,2,...,d\}^\mathbb{N}=\Omega$ and consider the compact metric space with metric $d$ where for
$x=(x_0,x_1,x_2,..)\in \Omega$ and $ y=(y_0,y_1,y_2,..)\in \Omega$
$$ d(x,y) = 2^{-N},$$
where $N$ is the smallest natural number $j\geq 0$, such that, $x_j\neq y_j$.

We also consider
$\{1,2,...,d\}^\mathbb{Z}=\hat{\Omega}$ and elements in $\hat{\Omega}$  are denoted by
$x=
(..., x_{-n}, ... ,x_{-1}\,|\,x_{0},x_1,..,x_n,..\,)$.

We will use the notation
$\overleftarrow{\Omega}= \{1,2,...,d\}^\mathbb{N}$ and
$\overrightarrow{\Omega}= \{1,2,...,d\}^\mathbb{N}$.

Given $x=
(..., x_{-n}, ... ,x_{-1}\,|\,x_{0},x_1,..,x_n,..\,)\in \hat{\Omega}$, we call $(..., x_{-n}, ... ,x_{-1}\,)\in  \overleftarrow{\Omega}$ the past of $x$ and  $(x_{0},x_1,..,x_n,..)\in \overrightarrow{\Omega}$ the future of $x$.

In this way we express  $\hat{\Omega}=\overleftarrow{\Omega}\times |\, \overrightarrow{\Omega} .$

Sometimes we denote
$$ (..., a_{-n}, ... ,a_{-1}\,|\,b_{0},b_1,..,b_n,..\,)=<a\,|\,b>,$$
where
$a= (..., a_{-n}, ... ,a_{-1})\in \overleftarrow{\Omega}$ and  $b=(\,b_{0},b_1,..,b_n,..)\in \overrightarrow{\Omega}$.

On $\hat{\Omega}$ we consider the usual metric $d$, in such way that for $x,y\in \hat{\Omega}$ we set
$$ d(x,y) = 2^{-N},$$
$N\geq 0$, where for
$$x=
(..., x_{-n}, ... ,x_{-1}\,|\,x_{0},x_1,..,x_n,..\,)\,\,\,,y=
(..., y_{-n}, ... ,y_{-1}\,|\,y_{0},y_1,..,y_n,..\,),$$ we have $x_j=y_j$, for all $j$, such that, $-N+1\leq j\leq N-1$
and, moreover $x_N\neq y_N$, or $x_{-N}\neq y_{-N}$.

The shift $\hat{\sigma}$ on $\hat{\Omega}=\{1,2,...,d\}^\mathbb{Z} $ is such that
$$\hat{\sigma} (..., y_{-n}, ... ,y_{-2},y_{-1}\,|\, y_{0},y_{1},....,y_{n},...)= (..., y_{-n}, ... ,y_{-2},y_{-1}, y_{0}\,|\,y_{1},....,y_{n},...).$$

On the other hand the shift $\sigma$ on $\Omega=\{1,2,...,d\}^\mathbb{N} $ is such that
$$\sigma ( y_{0},y_{1},....,y_{n},...)= (\,y_{1},....,y_{n},...).$$

\medskip

 A general equivalence relation $R$ on a space $X$  define classes and we will denote by $x \sim y$ when  two elements $x$ and $y$ are on the same class.
We denote by $[y]$ the class of $y\in X$.

\medskip

\begin{definition} Given an equivalence relation $\sim$ on $X $, where $X$ is any of the sets $\Omega, \,\hat{\Omega}$,  $(S^1)^\mathbb{N},$ or $S^1 \times S^1$, we denote by $G$ the subset
of $ X \times X  $,  containing all pairs $(x,y)$, where $x \sim y$. We call $G$ the groupoid associated to the equivalence relation $\sim$.

We also denote by $G^0$ the set  $\{(x,x)\,|\, x \in X\}\sim X $, where $X$ denote any of the sets  $\Omega, \,\hat{\Omega}$, $(S^1)^\mathbb{N},$ or $S^1 \times S^1$.

\end{definition}

{\bf Remark:}
There is a general definition of groupoid (see \cite{Con}) which assumes more structure but we will not need this here. For all  results we will consider  there is no need for an additional algebraic structure (on the class of each point). In this way we can consider a simplified definition of groupoid as it is above. Our intention is to study
$C^*$-algebras and Haar systems as a topic on measure theory (intersected with ergodic theory) avoiding
questions of algebraic nature.

\medskip

There is a future issue about the topology we will consider induced on $G$. One possibility is the product topology, which we call the standard structure, or,  a more complex one which   will be defined later on section \ref{C} (specially appropriate for some $C^*$-algebras).

\medskip

We will present several examples of dynamically defined groupoids.  The equivalence relation of most of our examples
is proper (see definition \ref{pro}).

\medskip

\begin{example} \label{exi}

For example consider on $\{1,2,...,d\}^\mathbb{N}$ the equivalence relation $R$ such that $x \sim y$,
if $x_j=y_j$, for all $j\geq 2$, when $x=(x_1,x_2,x_3,...)$ and  $y=(y_1,y_2,y_3,...)$. This defines a groupoid $G$.
In this case $G^0 =\Omega = \{1,2,...,d\}^\mathbb{N}$.

For a fixed $x=(x_1,x_2,x_3,...)$ the equivalence  class associated to $x$ is the set $ \{\,(j,x_2,x_3,...)$, $j=1,2..,d\,\}$.
We call this relation the {\bf bigger than two relation}.

\end{example}

\medskip

\begin{example}  \label{inv1}

Consider an equivalence relation $R$ which defines a  partition $\eta_0$ of $\{1,2,...,d\}^{\mathbb{Z}-\{0\}}=\hat{\Omega} $ such its elements are of the form
$$a\times \,|\, \overrightarrow{\Omega}=a  \times \{1,2,...,d\}^\mathbb{N}=  (..., a_{-n}, ... ,a_{-2},a_{-1}) \times \,|\, \{1,2,...,d\}^\mathbb{N},$$ where $a   \in \{1,2,...,d\}^\mathbb{N}= \overleftarrow{\Omega}$. This defines an equivalence relation $\sim$.

In this way two elements $x$ and $y$ are related if they have the same past.

There exists a bijection of classes of $\eta_0$ and points in $\overleftarrow{\Omega}$.

Denote $\pi=\pi_2:\hat{\Omega} \to \overleftarrow{\Omega}$ the transformation such that takes a point and gives as the result its class.

In this sense
$$\pi^{-1} (x) = \pi^{-1} ( (..., x_{-n}, ... ,x_{-1}\,|\,x_1,..,x_n,..)\,) ) = $$
 $$(..., x_{-n}, ... ,x_{-2},x_{-1}) \times \,|\, \overrightarrow{\Omega} \cong \Omega.$$

The groupoid obtained by this equivalence relation can be expressed as  $G=\{(x,y), \pi (x)=\pi(y)\}.$ In this way $x \sim y$ if they have the same past.

In this case the number of elements in each fiber is not finite.


Using the notation of page 46 of \cite{Con} we have $Y\subset \hat{\Omega}\times \hat{\Omega}$ and $X=\overleftarrow{\Omega}$.

In this case each class is associated to certain $a= (a_{-1},a_{-2},..) \in  \overleftarrow{\Omega}  =\{1,2,...,d\}^\mathbb{N}.$.

We use the notation $(a\,| \,x) $ for points on a class of the form
$$(a\,| \,x) =(..., a_{-n}, ... ,a_{-2},a_{-1}\,|\,x_{1},....,x_{n},...)\,\,.$$

\end{example}

\begin{example} \label{win} A particulary important equivalence relation $R$ on $\hat{\Omega}= \{1,2,...,d\}^{\mathbb{Z}}$ is the following: we say  $x \sim y$ if
$$x =(..., x_{-n}, ... ,x_{-2},x_{-1}\,|\, x_{0},x_{1},....,x_{n},...),\,\,$$
and,
$$y =(..., y_{-n}, ... ,y_{-2},y_{-1}\,|\, y_{0},y_{1},....,y_{n},...)\,\,$$ are such that there exists
$k\in \mathbb{Z},$ such that, $ x_j=y_j,$ for all $\,j\leq k  .$

The groupoid $G_u$ is defined by this relation $x \sim y$.

By definition the unstable set of the point $x\in\hat{\Omega}$ is the set
$$W^u (x)=\{y \in \hat{\Omega}\,,\,\,\text{such that}\,\,\, \lim_{n\to \infty} d(\, \hat{\sigma}^{-n}(x), \hat{\sigma}^{-n}(y)\,)\,=0\,\}.$$

One can show that the unstable manifold of $x \in \hat{\Omega}$ is the set
$$W^u (x)=\{y =(..., y_{-n}, ... ,y_{-2},y_{-1}\,|\, y_{0},y_{1},....,y_{n},...)\,|\,\text{there exists}\,$$
 $$ k\in \mathbb{Z},\,\, \text{such that}\,\, x_j=y_j, \text{for all}\,\,j\leq k  \}.$$

If we  denote by $G_u$ the groupoid   defined by the above relation, then, $x \sim y$, if and only if, $y \in W^u (x).$

An equivalence relation of this sort - for hyperbolic diffeomorphism -  was considered on \cite{Seg} and \cite{Man}.


\end{example}

\medskip

\begin{example} \label{win1} An equivalence relation  on $
 \overrightarrow{\Omega}= \{1,2,...,d\}^\mathbb{N}$ similar to the previous one is the following: we say  $x \sim y$ if
$$x =(\, x_{0},x_{1},....,x_{n},...),\,\,$$
and,
$$y =(y_{0},y_{1},....,y_{n},...)\,\,$$ are such that there exists
$k\in \mathbb{N},$ such that, $ x_j=y_j,$ for all $\,j\geq k  .$

\end{example}

\begin{example} \label{win590} Another equivalence relation  on $
 \overrightarrow{\Omega}$ is the following: fix $k \in \mathbb{N}$,  and we say  $x \sim_k y$, if when
$$x =(\, x_{0},x_{1},....,x_{n},...),\,\,$$
and,
$$y =(y_{0},y_{1},....,y_{n},...)\,\,$$ we have $ x_j=y_j,$ for all $\,j\geq k  .$

In this case each class has $d^k$ elements.

\end{example}

\begin{example} \label{homore}

Given $x,y \in \hat{\Omega}= \{1,2,...,d\}^\mathbb{Z} $, we say that $x\sim y$ if
$$ \lim_{k \to + \infty} d(\hat{\sigma}^k x, \hat{\sigma}^k y)=0$$
$$\text{and}$$
\begin{equation} \label{equivalencee}
\lim_{k \to - \infty} d(\hat{\sigma}^k x, \hat{\sigma}^k y)=0.
\end{equation}

This means there exists an $M\geq 0$, such that, $x_j=y_j$ for $j> M$, and, $j< -M$. In other words, there are only a finite number of $i$'s such that $x_{i} \neq y_{i}$. This is the same to say that $x$ and $y$ are homoclinic.

For example in $\hat{\Omega}=\{1,2\}^\mathbb{Z}$ take
$$ x= (...,x_{-n},...,x_{-7},1,2,2,1,2,2\,|\, 1,2,1,2,1,1,x_7,...x_{n},..)$$
and
$$ y= (...,y_{-n},...,y_{-7},1,2,2,1,2,2\,|\, 1,2,1,1,1,2,y_7,...y_{n},..)$$
where $x_j=y_j$ for $|j| \geq 7.$

In this case $x \sim y$.

This relation is called the {\bf homoclinic relation} on $\hat{\Omega}.$  It was considered for instance by D. Ruelle and N. Haydn in \cite{Ru1} and \cite{HD} for hyperbolic diffeomorphisms and also on more general contexts (see also \cite{LM2}, \cite{Meye} and \cite{Bis}  for the symbolic case). \end{example}

\begin{example} \label{baker} Consider an expanding transformation $T:S^1 \to S^1$, of degree two,  such that $\log T\,'$ is Holder
and $\log T\,'(a) > \log \lambda >0$, $a \in S^1$, for some $\lambda >1$.

Suppose $T(x_0)=1,$ where $0<x_0<1$. We say that $(0,x_0)$ and $(x_0,1)$ are the domains of injectivity of $T$.

Denote $\psi_1:[0,1) \to [0,x_0)$ the first inverse branch of $T$ and $\psi_2:[0,1) \to[x_0,1]$ the second inverse branch of $T$.

In this case for all $y$ we have $T \circ \psi_1(y)=y$ and  $T \circ \psi_2(y)=y$.

The {\bf associated  $T$-Baker map} is the transformation $F:S^1 \times S^1$ such that  satisfies for all $a,b$ the following rule:

1) if $0\leq b < x_0$
$$ F(a,b)= (\psi_1(a), T(b)),$$

and

2) if $x_0 \leq b < 1$
$$ F(a,b)= (\psi_2(a)\, ,T(b)).$$

In this case we take as partition the one associated to (local) unstable manifolds for $F$, that is, sets of the form $W_a = \{\,(a,b)\,|\, b \in S^1\} $, where $a \in S^1$.

Given two points $z_1,z_2\in S^1 \times S^1$ we say that they are related if the  the first coordinate is equal.

On $S^1\times S^1$ we use the distance $d$ which is the product of the usual arc length distance on $S^1$.

The bijection $F$ expands vertical lines and contract horizontal lines.

As an example one can take $T(a)= 2 a \, $ (mod 1) and we get (the inverse of)  the classical Baker map (see \cite{Bak}).

One can say that the dynamics of such $F$ in some sense looks like the one of an Anosov diffeomorphism.

\end{example}

\begin{example} \label{XY} The so called generalized $XY$ model consider space
$(S^1)^\mathbb{N}$, where $S^1$ is the unitary circle and the shift acting on it (see ).

We can consider the equivalence relation $R$ such that $x \sim y$,
if $x_j=y_j$, for all $j\geq 2$, when $x=(x_1,x_2,x_3,...)$ and  $y=(y_1,y_2,y_3,...)$. This defines a groupoid $G$.
In this case $G^0 =(S^1)^\mathbb{N}$.

For a fixed $x=(x_1,x_2,x_3,...)$ the equivalence  class associated to $x$ is the set $ \{\,(a,x_2,x_3,...)$, $a\in S^1\,\}$.
We call this relation the {\bf bigger than two relation for the $XY$ model} and $G$ the {\bf standard $XY$ groupoid} over $(S^1)^\mathbb{N}$.

\end{example}

\section{Kernels and transverse functions} \label{ktran}

\medskip

A general reference for the material  of this section is \cite{Con} (see also \cite{Kas}).

We consider over $G\subset X\times X$ the  Borel sigma-algebra $\mathcal{B}$ (on $G$) induced by the natural
product topology on  $X \times X$ and the metric $d$ on $X$ (\cite{FM1} and \cite{FM2} also consider this sigma algebra). This will be fine for the setting of von Neumann algebras.
Later, another sigma-algebra will be considered for the setting $C^*$-algebras.

We point out  that the only sets $X$ which we are interested are of the form $\hat{\Omega}$,  $\Omega$,  $(S^1)^\mathbb{N},$ or $S^1 \times S^1$.

We denote
$\mathcal{F}^+( G) $ the space of Borel measurable functions $f: G \to [0,\infty)$ (a function of two variables $(a,b)$).

$\mathcal{F}( G) $ is  the space of Borel measurable functions $f: G \to \mathbb{R}$. Note that $f(x,y)$ just make sense if $x \sim y$.

There is a natural  involution on $\mathcal{F}( G) $ which is $f \to \tilde{f}$, where $\tilde{f}(x,y)= f(y,x).$

 We also denote
$\mathcal{F}^+( G^0) $ the space of Borel measurable functions $f: G^0 \to [0,\infty)$ (a function of  one variable $a$).

There is a natural identification of functions $f: G^0 \to \mathbb{R}$, of the form $f(x)$, with functions $g: G \to \mathbb{R}$ which depend only on the first coordinate, that is $g(x,y)=f(x)$. This will be used without mention, but if necessary we write $(f\circ P_1)(x,y)=f(x)$ and $(f\circ P_2)(x,y)=f(y)$.

\begin{definition} A {\bf measurable groupoid} $G$ is a groupoid with the topology induced by the  product topology over $X\times X$, such that, the following functions are measurable for the Borel sigma-algebra:

$P_1(x,y)=x$, $\,P_2(x,y)=y$,  $\,h(x,y)=(y,x)$ and $Z (\,(x,s),(s,y)\,)= (x,y),$ where $Z:\, \{\,(\,(x,s),(r,y) \,)\,\,|\,\, r=s\,\}\,\subset G \times G\,\to \,G.$

\end{definition}

Now, we will present the definition of kernel (see beginning of section 2 in \cite{Con}).

\begin{definition}
A {\bf G-kernel} $\nu$ on the measurable groupoid $G$ is an application of $G^0$ in the space of measures over the sigma-algebra $\mathcal{B}$, such
that,

1) for any $y \in G^0$, we have that $\nu^y$ has support on $[y]$,

and

2) for any $A \in \mathcal{B}$, the function $y \to \nu^y (A)$ is measurable.

 The set of all $G$- kernels is denoted by $\mathcal{K}^+.$

\end{definition}

\begin{example} \label{btt}
As an example consider for the case
of the  groupoid $G$ associated to the  bigger than two relation,
the measure $\nu^y$, for each $y=(y_1,y_2,y_3,...)$, such that $ \nu (j,y_2,y_3,...)= 1$, $j=1,2,...,d$.
In other words we are using the counting measure on each class.We call this the standard $G$-kernel for the the  bigger than two relation.

More precisely, the {\bf counting measure} is such that $\nu^y (A) = \# (A \cap [y]),$ for any $A \in \mathcal{B}.$

\end{example}

\begin{example} \label{bttd}
Another possibility is to consider the $G$-kernel such that
$\nu^y$, for each $y=(y_1,y_2,y_3,...)$, is such that $ \nu (j,y_2,y_3,...)= \frac{1}{d}$,
We call this the normalized standard $G$-kernel for the the  bigger than two relation.
\end{example}
\medskip

\begin{example} \label{del}
Given any groupoid $G$ another example of kernel is the {\bf delta kernel} $\nu$ which is the one such that for any $y \in G^0$ we have that
$\nu^y (dx) = \delta_y (d\,x)$, where $\delta_y$ is the delta Dirac on $y$. We denote by $\mathfrak{d}$ such kernel.
\end{example}

\medskip

We denote by  $\mathcal{F}_\nu (G)$ the set of $\nu$-integrable functions.

\medskip

\begin{definition}\label{intkernel}

Given a $G$-kernel $\nu$ and
an integrable function $f\in \mathcal{F}_\nu (G)$ we denote by  $\nu  ( f)$ the function in $\mathcal{F} (G^0)$ defined by
$$ \nu  (f)\,\,(y) =\int f\,(s,y)\, \nu^y (ds), \,\,\, y \in G^0.$$

\end{definition}

A kernel $\nu$ is characterized by the law
$$f\in \mathcal{F}_\nu (G)\,\,\,\to\,\,\,\nu  ( f)\in \mathcal{F} (G^0).$$
\medskip

In other words, for a kernel $\nu$ we get
$$ \nu : \,  \mathcal{F}_{\nu} (G) \to  \mathcal{F} (G^0).$$

By notation given a kernel $\nu$ and a positive $f\in \mathcal{F}_\nu (G)$ then the kernel $ f\, \nu$ is the one defined by $f(x,y) \nu^y (dx).$ In other words the action of the kernel $ f\, \nu$ get rid of the first coordinate:

$$ h(x,y) \,\,\to\,\, \int h(s,y) f\,(s,y)\, \nu^y (ds).$$

In this way if $f\in \mathcal{F}_\nu (G^0)$ we get $f(x) \nu^y(dx)$.
\medskip

Note that $\nu(f)$ is a function and $f\, \nu$ is a kernel.

\medskip

\begin{definition} \label{wi2}

 A {\bf transverse function} is a $G$-kernel $\nu$, such that, if $x \sim y$, then, the finite measures  $\nu^y$ and $\nu^x$ are the same. The set of transverse functions for $G$ is denoted by $\mathcal{E}^+.$ We call probabilistic transverse function any one such that for each $y \in G^0$ we get that $\nu^y$  is a probability on the class of $y$.






\medskip
\end{definition}

The above means that
$$ \int f(a)  \nu^x (d\,a)= \int f(a)  \nu^y (d\, a),$$
if $x$ and $y$ are related. In the above we have $x \sim y \sim a.$

\medskip

\begin{remark}
The above equality implies that a transverse function is left (and right) invariant. Together with the conditions defining a $G$-kernel, we have that $\nu$ is {\bf a Haar system} (see \cite{Ren0}) in a {\bf measurable sense}. In what follows we use {\bf transverse function and Haar system as synonyms}.
\end{remark}


The standard $G$-kernel for the  bigger than two relation (see example \ref{btt})  is a transverse function.

The normalized standard $G$-kernel for the the  bigger than two relation (see example \ref{bttd})  is a probabilistic transverse function.

If we consider the equivalence relation such that each point is related just to itself, then the transverse functions
can be identified with the positive functions defined on $X$.

The difference between a function and a transverse function is that the former takes values on the set of real numbers and the later on the set of measures.

\bigskip

If $\nu$ is transverse, then $\nu^x=\nu^y$ when $x \sim y$, and we have from definition \ref{kern*func}:

\begin{equation}\label{eq1}
(\nu  * f)(x,y)= \int f\,(x,s)\, \nu^x (ds) = \nu (\tilde f)(x),\,\,\, \forall y \sim x
\end{equation}
and,
\begin{equation}\label{eq2}
(f * \nu)(x,y)= \int f\,(s,y)\, \nu^y (ds) = \nu(f)(y) ,\,\,\, \forall x \sim y.
\end{equation}

\medskip

\begin{definition} The pair $(G,\,\nu)$, where $\nu\in\mathcal{E}^+$, is called the {\bf measured groupoid for the transverse function $\nu$}. We assume  any $\nu$ we consider is such that $\nu^y$ is not the zero measure for any $y$.
\end{definition}

\medskip

In the case $\nu$ is such that,  $\int \nu^{y} (ds)=1$, for any $y \in G^0$, the Haar system will be called a probabilistic Haar system.

\medskip

Note that the delta kernel $\mathfrak{d}$ is not a transverse function.

\medskip

\medskip

 Given a measured groupoid $(G,\nu)$ and two measurable functions  $f,g \in  \mathcal{F}_\nu (G)$, we  define $(f \,  \underset{\nu}{*}\, g)=h$ in such way that for any $(x,y)\in G$
$$(f \,  \underset{\nu}{*}\, g)(x,y)=
\int \,g(x,s ) \,f (s,y) \,\,\nu^{y} (ds) =h(x,y).$$

$(f \,  \underset{\nu}{*}\, g)$ is called the {\bf  convolution} of the functions $f,g$ for the measured groupoid $(G,\nu)$.


\begin{example} \label{exi25} Consider the groupoid $G$ of Example \ref{exi} and the family $\nu^y$, $y \in \{1,2,...,d\}^\mathbb{N}$, of measures (where each measure $\nu^y$ has support on the equivalence class of
$y$), such that, $\nu^y$ is the counting measure. This defines a transverse function (Haar system) called the {\bf standard Haar system}.
\end{example}

\begin{example} \label{exi26} Consider the groupoid $G$ of Example \ref{exi} and the normalized standard family $\nu^y$, $y \in \{1,2,...,d\}^\mathbb{N}$. This defines a transverse function called the {\bf normalized standard Haar system}.

More precisely the family $\nu^y$, $y \in \{1,2,...,d\}^\mathbb{N}$, $y=(y_1 ,y_2,y_3,...)$, of probabilities on the set
$$\{\,(a,y_2,y_3,...), a \in \{1,2,..,d\}\,\,\},$$
is  such that, $\nu^y(\, \{(a ,y_2,y_3,...)\,\})= \frac{1}{d}$, $a \in \{1,2,..,d\}$
\end{example}

\begin{example} \label{win591} In example \ref{win590} in which $k$ is fixed consider the transverse
function $\nu$ such that for each $y \in G^0$, we get that  $\nu^y$ is the counting measure on the set of points $x\sim_k y$.

\end{example}

\begin{example} \label{exi2} Suppose $J:\{1,2,...,d\}^\mathbb{N}\to  \mathbb{R}$ is continuous positive function such that for any $x\in \Omega$ we have that
$\sum_{a=1}^d J(ax)=1$. For the groupoid $G$ of Example \ref{exi}, the family $\nu^y$, $y \in \{1,2,...,d\}^\mathbb{N}$, of probabilities on
$\{\,(a,y_2,y_3,...), a \in \{1,2,..,d\}\,\,\}$, such that, $\nu^y( a ,y_2,y_3,...)= J(a,y_2,y_3,...)$, $a \in \{1,2,..,d\}$ defines a Haar system. We call it the {\bf probability Haar system associated to $J$}.

Example \ref{exi26} is a particular case of the present example.

\end{example}

\begin{example} \label{inv2} On the groupoid over $\{1,2,...,d\}^\mathbb{Z}$ described on example \ref{inv1}, where we consider the notation: for each class specified by $a\in \overleftarrow{\Omega}$ the general element in the class is given by
$$(a\,| \,x) =(..., a_{-n}, ... ,a_{-2},a_{-1}\,|\,x_{1},....,x_{n},...),$$
where $x \in \overrightarrow{\Omega}$.

Consider a fixed probability $\mu $ on $\overrightarrow{\Omega}.$ We define the transverse function $\nu^a (dx) = \mu(dx)$ independent of $a$.

\end{example}

\begin{example} \label{baba} In the example \ref{baker} we consider the partition of $S^1 \times S^1$ given by the sets
$W_a = \{\,(a,b)\,|\, b \in S^1\} $, where $a \in S^1$.  For each $a\in S^1$, consider a probability $\nu^a (d\,b)$ over $\{\,(a,b)\,|\, b \in S^1\} $ such that for any Borel set $ K \subset S^1 \times S^1$ we have that $a \to \nu^a (K)$ is measurable. This defines a probabilistic transverse function and a Haar system.

Consider a continuous function $A: S^1 \times S^1 \to \mathbb{R} $. For each $a$ consider the kernel
$\nu^a $ such that $\int f(b) \nu^a (db) = \int f(b) \, e^{A(a,b)} db$, where $d b$ is the Lebesgue measure. This defines a  transverse function.

We call the {\bf standard Haar system on $\mathbf{S^1 \times S^1}$ } the case where for each $a$ we consider as the probability $\nu^a (d\,b)$ over $\{\,(a,b)\,|\, b \in S^1\} $ the Lebesgue probability on $S^1$.

\end{example}

We will present several properties of kernels and transverse functions on Section \ref{non}.

\bigskip

A question of notation: for a fixed groupoid $G$ we will  describe now for the reader the common  terminology on the field (see \cite{Con}, \cite{Kas}, \cite{Ren0} and \cite{Ren2}). It is usual to denote a general pair $(x,y)\in G$ by $\gamma$ (of related elements $x$, $y$). The $\gamma$ is called the directed arrow from $x$ to $y$. In this case we call $s(\gamma) =x$ and $r(\gamma)=y$ (see \cite{Ou} for a more detailed description of the arrow's setting).

Here, for each pair of related elements $(x,y)$ there exist an unique directed arrow $\gamma$  satisfying $s(\gamma) =x$ and $r(\gamma)=y$. Note that, since se are dealing with equivalence relations, $(y,x)$ denotes another arrow. In category language: there is a unique morphism $\gamma$ that takes $\{x\}$ to $\{y\}$, whenever $x$ and $y$ are related, and this morphism is associated in a unique way to the pair $(x,y)$.

In this notation $r^{-1} (y)$ is the set of all arrows that end in $y$. This is in a bijection with all elements on the same class of equivalence of $y$. We call  $r^{-1} (y)$ the fiber over $y$. If $x \sim y$, then $r^{-1} (y)=r^{-1} (x)$.

We adapt the notation in \cite{Con} and \cite{Kas} to our notation. We use here the expression $(s,y)$ instead of $\gamma\, \gamma'$. This makes sense considering that $\gamma=(x,y)$ and $\gamma'=(s,x)$. We use the expression $(y,s)$ for $(\gamma')^{-1}\, \gamma$, where in this case, $\gamma=(x,y)$ and $\gamma'=(s,y)$, and, finally,  $\nu^y ( \gamma')$ means $ \nu^y (ds)$ for $\gamma'=(s,y)$.

In the case of the  groupoid $G$ associated to the  bigger than two relation we have
for each $x=(x_1,x_2,x_3,...)$  the property $ r^{-1}(x)=\{\,(j,x_2,x_3,...)$, $j=1,2..,d\,\}$.

The terminology of arrows will not be essentially used here. It was introduced just for the reader to make a parallel (a dictionary)  with the one commonly used on papers on the topic.

Using the terminology of arrows Definition \ref{wi2} is equivalent to say that: if, $\gamma=(x,y)=(s(\gamma),r(\gamma))$, then,
$$ \nu^{ y} = \gamma\, \nu^{x}.$$

\medskip

\section{Quasi-invariant probabilities} \label{quasi}

\begin{definition} A function $\delta:G \to \mathbb{R}$ such that
$$ \delta (x,z) = \delta (x,y)\, \delta (y,z)  ,$$
for any $(x,y)$, $(y,z)\in G$ is called a {\bf modular function} (also called a multiplicative {\bf cocycle}).

In the arrow notation this is equivalent to say that
$$ \delta (\gamma_1 \gamma_2) = \delta (\gamma_1) \delta (\gamma_2). $$

\end{definition}

\bigskip

Note that $\delta(x,y)\, \delta(y,y)= \delta(x,y)$ and it follows that for any $y$ we have $\delta(y,y)=1.$ Moreover,
$\delta(x,y) \,\delta(y,x)= \delta(x,x)=1$ is true. Therefore, we  get $\tilde{\delta}=\delta^{-1}.$

\medskip

\begin{example} \label{bib1} Given $W: G^0 \to \mathbb{R}$, $W(x)> 0, \forall x$, a natural way to get a modular function is to consider
$\delta(x,y)= \frac{W(x)}{W(y)}.$ In this case we say that the modular function  is derived from $W$.
\end{example}

\begin{example} \label{baker1} In the case of example \ref{baker} the equivalence relation is:
given two points $z_1,z_2\in S^1 \times S^1$ they are related if the first coordinate is equal.

Consider a expanding transformation $T$ and the associated Baker map $F$. Note que $F^n (a,b)= (*, T^n (b))$ for some point $*$.

Given two points $z_1\sim z_2$, for each $n$
there exist  $z_1^n$ and   $z_2^n$, such that, respectively, $F^n(z_1^n)=z_1$ and  $F^n(z_2^n)=z_2$, and $z_1^n \sim z_2^n$.

For each pair $z_1=(a,b_1)$ and $z_2=(a,\,b_2)$, and $n \geq 0$,  the elements  $z_1^n, z_2^n$ are of the form
$z_1^n= (a^n, b_1^n),$ $ z_2^n=(a^n, b_2^n).$

In this case $ T^n (b_1^n)= b_1$ and  $ T^n (b_2^n)= b_2.$

Note also that $T^n (a)= a^n$.

The distances between $b_1^n $ and $b_2^n$ are exponentially decreasing with $n$.

We denote
$$ \delta(z_1,z_2)=\Pi_{j=1}^\infty\,\, \,\frac{T\,' ( b_1^n )}{T\, ' (b_2^n )}< \infty .$$

This product is well defined because
$$\sum_n \, \log  \frac{T\,' ( b_1^n )}{T\, ' (b_2^n )} =
\sum_n \,[ \log  T\,' ( b_1^n )  \,-\, \log T\, ' (b_2^n ) \,]$$
converges. This is so because $\log  T\,'$ is Holder and for all $n$ we have   $|b_1^n - b_2^n|< \lambda^{-n},$ where $T' (x) > \lambda>1$ for all $x$.

This $\delta$ is a cocycle.

In the case of example \ref{baba}
considered a Holder continuous function $A(a,b)$, where $A: S^1 \times S^1 \to \mathbb{R} $.

Define for $z_1=(a,b_1)$ and $z_2=(a,\,b_2)$
$$ \delta(z_1,z_2)=\Pi_{j=1}^\infty\,\, \,\frac{e^{A ( z_1^n )}}{e^{A(z_2^n )}}.$$

The modular function $ \delta(z_1,z_2)$ is well defined because $A$ is Holder.

We will show that $\delta$ can be expressed in the form of example \ref{bib1}. Indeed, fix a certain $b_0\in S^1$, then, taking $z_1=(a,b_1) $ consider $z_0=(a,b_0) $. We denote  in an analogous way  $z_1^n$ and $z_0^n$ the ones such that  $F^n(z_1^n)=z_1$ and  $F^n(z_0^n)=z_0$.

Define $V: G^0 \to \mathbb{R}$ by
\begin{equation} \label{Lar} V(z_1)=\Pi_{j=1}^\infty\,\, \,\frac{e^{A ( z_1^n )}}{e^{A(z_0^n )}}.
\end{equation}

$V$ is well defined and if $z_1\sim z_2$ we get that
$$ \delta(z_1,z_2)= \frac{V(z_1)}{V(z_2)}.$$

We will show later (see Proposition \ref{rod}) that $V(a,b)$ does not depend on $a$, and then we can write $V(b)$, and finally
$$ \delta(z_1,z_2)= \frac{V(b_1)}{V(b_2)}.$$

\end{example}

\begin{example}
Consider a  fixed Holder function $\hat{A}: \{1,2,...,d\}^\mathbb{Z}\to \mathbb{R}$ and the groupoid given by the equivalence relation of example \ref{win}. Denote for any $(x,y)$
$$ \delta(x,y)=\Pi_{j=1}^\infty\,\, \,\frac{\hat{A} ( \hat{\sigma}^{-j} (s(\gamma)) )}{\hat{A} ( \hat{\sigma}^{-j} (r(\gamma)) )}= \Pi_{j=1}^\infty\,\, \,\frac{\hat{A} ( \hat{\sigma}^{-j} (x) )}{\hat{A} ( \hat{\sigma}^{-j} (y) )}  .$$

The modular function $\delta$ is well defined because $\hat{A}$ is Holder. Indeed, this follows from  the bounded distortion property.

In a similar way as in the last example one can show that such $\delta$ can be expressed on the form of example  \ref{bib1}.

\end{example}

\medskip

\begin{definition} \label{der} Given a measured groupoid $G$  for the transverse function $\nu$ we say that a probability $M$ on $G^0$  is {\bf quasi-invariant} for $\nu$ if there exist a modular function $\delta:G \to \mathbb{R}$, such that,
for any integrable function  $f:G \to \mathbb{R}$ we have
\begin{equation} \label{kwe}
 \int \, \int f(s,x) \nu^x (ds) d M(x)= \int \, \int f(x,s) \,\delta^{-1}(x,s)\,\nu^x (ds) d M(x).
\end{equation}
\end{definition}

In a more accurate way we say that $M$ is quasi-invariant for the transverse function
$\nu$ and the modular function $\delta$.


\medskip

Note that if $\delta(x,s) = \frac{B(x)}{B(s)} $  we get that the above condition (\ref{kwe}) can be written as
\begin{equation} \label{kwe1} \int \, \int f(s,x)  B(s) \nu^x (ds) d M(x)= \int \, \int f(x,s) \,B(s)\,\nu^x (ds) d M(x).
\end{equation}

Indeed, in (\ref{kwe}) replace $f(s,x)$ by $ B(s) f(s,x)$.

\medskip

Quasi-invariant probabilities will be also described as the ones which satisfies the
so called the KMS condition (on the setting of  von Neumann algebras, or $C^*$-algebras) as we will see later on section \ref{von1}.

As an extreme example consider the equivalence relation such that each point is related to just itself. In this case
a modular function $\delta$ takes only the value $1$. Given any transverse function $\nu$ the condition

\begin{equation} \label{kwe99}
 \int \, \int f(s,x) \nu^x (ds) d M(x)= \int \, \int f(x,s) \,\delta^{-1}(x,s)\,\nu^x (ds) d M(x)
\end{equation}
is satisfied by any probability $M$ on $X$. In this case the set of probabilities is the set of  quasi-invariant probabilities.

\medskip

\begin{example} \label{SBR}  {\bf Quasi invariant probability  and the SBR probability for the Baker map }

\medskip

We will present a particular example where we will compare the probability $M$ satisfying the quasi invariant condition with  the so called SBR probability. We will consider a different setting of the case  described on \cite{Seg}  (considering Anosov systems) which, as far as we know, was never published.

We will show that the {\bf  quasi invariant probability is not the SBR probability}.

We will address later on the end of this example the kind of questions discussed on  \cite{Seg} and \cite{Man}.

We will consider the groupoid of example \ref{baker}, that is, we consider the equivalence relation:
given two points $z_1,z_2\in S^1 \times S^1$ they are related if the first coordinate is equal.

In example \ref{baker} we consider an expanding transformation $T:S^1 \to S^1$ and
$F$ denotes the associated $T$-Baker map. The associated SBR probability is the only absolutely continuous $F$-invariant probability over
$S^1 \times S^1$.

The dynamical action of $F$ in some sense looks like the one of an Anosov diffeomorphism.


Consider the measured groupoid $(G\,,\nu)$ where in each vertical fiber over the point $a$ we set $\nu^a$  as the Lebesgue probability
$ db$ over the class $(a,b), $ $ 0\leq b\leq 1.$

This groupoid corresponds to the local unstable  foliation for the transformation $F$.

We fix a certain point $b_0\in(0,1)$.
For each pair $x=(a,b)$ and $y=(a,\,b_0)$, where $a,b \in S^1$, and $n \geq 0$, the elements  $z_1^n, z_2^n$, $n \in \mathbb{N}$, are such that $F^n (z_1^n)=x=(a,b)$ and  $F^n (z_2^n)=y=(a,b_0)$. Note that they are of the form
$z_1^n= (a^n, b^n),$ $ z_2^n=(a^n, s^n)$.
We use the notation $z_1^n(x)$, $b^n( x)$, $n \in \mathbb{N}$, to express the dependence on $x$.

We denote for $x \in S^1 \times S^1$
$$V(x)= V(a,b)=\Pi_{n=1}^\infty\,\, \,\frac{T' ( b^n (x))}{T' (s^n) }=   \Pi_{n=1}^\infty\,\, \,\frac{T' ( b^n (a,b))}{T' (s^n) }< \infty .$$

This is finite because  $s^n$ and
$b^n (x)$ are on the same domain of injectivity of $T$ for all $n$ and $T'$ is of Holder class.

In a similar fashion as in \cite{Seg} we define $\delta $ by the expression
$$ \delta(\,(a,y_1)\,,\, (a,y_2)\,)=\frac{V(a,y_1)}{V(a,y_2)}=\frac{V(y_1)}{V(y_2)}= \Pi_{n=1}^\infty\,\, \,\frac{T' ( b^n (a,y_1))}{  T' ( b^n (a,y_2))},$$
where $(a,y_1)\,\sim\, (a,y_2)$.

Consider the probability $M$ on $S^1\times S^1$ given by
$$dM(a,b)= \frac{V(a,b)}{ \int V(a,c) dc}\, db\, da.$$

The density $\psi(a,b) =  \frac{V(a,b)}{ \int V(a,c) dc}\,$
satisfies the equation
\begin{equation} \label{ali} \psi (a,b) \frac{1 }{T ' (b)} = \psi (F(a,b)).\end{equation}

Denote $F(a,b) = (\tilde{a}, \tilde{b})$, then, it is known that  the density $ \varphi(a,b)$ of the SBR probability for $F$ satisfies
\begin{equation} \label{ali1} \varphi (a,b) \frac{T ' (\tilde{a}) }{T ' (b)} = \varphi (F(a,b)).\end{equation}

This follows from the $F$-invariance of the SBR

Therefore, $M$ is not the SBR probability - by uniqueness of the SBR.

\medskip

We will show that $M$ satisfies the quasi invariant condition.

Note that
$$\int \, \int f(\,(a,b),(a,s)\,) \,\,\nu^a (d\,s) \,\,d M(a,b)= $$
$$  \int \,\int \, \int f(\,(a,b),(a,s)\,) \, \,  \,\, \frac{V(a,b)}{ \int V(a,c) dc}\,ds\, db\, da.$$

On the other hand

$$\int \, \int f(\,(a,s),(a,b)\,) \,\,\frac{ V(a,s)}{ V(a,b)\,}\,\,\nu^a (d\,s) \,\,d M(a,b)=
$$
$$\int \int  \, \int f(\,(a,s),(a,b)\,) \,\,\frac{ V(a,s)}{ V(a,b)\,}\, \, \,\, \frac{V(a,b)}{ \int V(a,c) dc}\,ds \, db\, da=
$$
$$\int \int \, \int f(\,(a,s),(a,b)\,) \,\, \,\,  \,\, \frac{V(a,s)}{ \int V(a,c) dc}\, ds \,db\, da.
$$

If we exchange the variables $b$ and $s$, and using Fubini's theorem, we get that $M$ satisfies the quasi invariant condition.
\medskip

The relation of quasi-invariant probabilities and transverse measures is described on section \ref{non}.

The result considered on Theorem 6.18 in \cite{Seg} for an Anosov diffeomorphism concerns transverse measures and cocycles. \cite{Seg} did not mention quasi-invariant probabilities.

Note that from equations (\ref{ali}) and (\ref{ali1}) one can get that the conditional disintegration along unstable leaves of both the SRB and the quasi-invariant probability $M$
are equal (see page 533 in \cite{LY}).

Using the relation of quasi-invariant probabilities, cocycles and transverse measures
one can say that one of the main claims in \cite{Seg} (see Theorem 6.18) and \cite{Man} (both considering the case of Anosov Systems) can be expressed in some sense  via the above mentioned property about conditional disintegration along unstable leaves (using the analogy with the case of the above Baker map $F$).

\end{example}

In section \ref{exaquasi} we will present more examples of quasi-stationary probabilities.

\section{von Neumann Algebras derived from measured groupoid}\label{von1}

We refer the reader  to \cite{AnaR}, \cite{Kas} and \cite{Con} as  general references for von Neumann algebras related to groupoids.

Here $X\sim G^0$ will be either $\hat{\Omega}$, $\Omega$ or $S^1 \times S^1$. We will denote by $G$ a general groupoid obtained by an equivalence relation $R$.

\begin{definition}\label{mult.algebra}
 Given a measured groupoid $G$ for the transverse function $\nu$  and two measurable functions  $f,g \in  \mathcal{F}_\nu (G)$, we  define the convolution  $(f \,  \underset{\nu}{*}\, g)=h$, in such way that, for any $(x,y)\in G$
$$(f \,  \underset{\nu}{*}\, g)(x,y)=
\int \,g(x,s ) \,f (s,y) \,\,\nu^{y} (ds) =h(x,y).$$

In the case there exists a multiplicative neutral element for the operation $*$ we denote it by $\mathfrak{1}.$
\end{definition}

The above expression in some sense resembles the way we get a matrix as the product of two matrices.

For a fixed Haar system $\nu$ the product  $\underset{\nu}{*}$ defines an algebra  on the vector
space of $\nu$-integrable functions  $\mathcal{F}_\nu (G)$.

As usual function of the form $f(x,x)$ are identified with functions $f: G^0 \to \mathbb{R}$ of the from $f(x).$

\begin{example} \label{mama}
In the particular case where $\nu^y $ is the counting measure on the fiber over $y$  then
$$(f \,  \underset{\nu}{*}\, g)(x,y)= \sum_s g(x,s)\,f(s,y) .$$

Denote by $I_\Delta$ the indicator function of the diagonal on $G^0 \times G^0$.
In this case, $I_\Delta$ is the neutral element for the product $\underset{\nu}{*}$ operation.

In this case $ \mathfrak{1}= I_\Delta$.

Note that
$I_\Delta$ is measurable but generally not continuous. This is fine for the von Neumann algebra setting.  However, we will
need a different topology (and $\sigma$-algebra) on $G^0 \times G^0$ - other than the product topology - when considering the unit $ \mathfrak{1}= I_\Delta$ for the $C^*$-algebra setting (see \cite{Dea}, \cite{Ren1.5}, \cite{Ren2}). This will be
more carefully explained on section \ref{C}.

\end{example}

{\bf Remark:} The indicator function of the diagonal on $G^0 \times G^0$ is not always the multiplicative neutral element on the von Neumann algebra obtained from a general Haar system $(G,\nu).$

\medskip

\begin{example} \label{exi4} Another example: consider the standard Haar system of example \ref{exi25}.

In this case

$$(f \,  \underset{\nu}{*}\, g)(x,y)=
\int \,g(x,s ) \,f (s,y) \nu^{y} (ds) =$$
$$ \frac{1}{d} \sum_{a=1}^d\, g(x, \,(\, a,x_2,x_3,...)\,) \,\,f(\,(a,x_2,x_3,...   )\,, y\,) \,\,  = h(x,y).$$

The neutral element is $d \,I_\Delta= \mathfrak{1}.$

\end{example}

\begin{example} \label{exi3} Suppose $J:\{1,2,...,d\}^\mathbb{N}\to  \mathbb{R}$ is a continuous positive function such that for any $x\in \Omega$ we have that
$\sum_{a=1}^d J(ax)=1$. The measured groupoid $(G,\nu)$ of Example \ref{exi2}, where  $\nu^y$, $y \in \{1,2,...,d\}^\mathbb{N}$, is such that given $f,g:G \to \mathbb{R}$, we have for any $(x,y)\in G$, $x=(x_1,x_2,x_3,...)$, $y=(y_1,x_2,x_3,..)$ that

$$(f \,  \underset{\nu}{*}\, g)(x,y)=
\int \,g(x,s ) \,f (s,y) \nu^{y} (ds) =$$
$$ \sum_{a=1}^d\, g(x, \,(\, a,x_2,x_3,...)\,) \,\,f(\,(a,x_2,x_3,...   )\,, y\,) \,\, J(a,x_2,x_3,...) = h(x,y).$$

Note that $x_j=y_j$ for $j\geq 2.$

Suppose that $f$ is such that for any string $(x_2,x_3,...)$ and $a\in \{1,2,.,d\}$ we get
$$ f(\,(a,x_2,x_3,...   )\,,(a,x_2,x_3,...   ) \,) \,=\,\frac{1}{ J(a,x_2,x_3,...)},$$
and, $a,b \in \{1,2,.,d\}$, $a \neq b$
$$ f(\,(a,x_2,x_3,...   )\,,(b,x_2,x_3,...   ) \,) \,=\,0.$$

In this case the neutral multiplicative element is $\mathfrak{1}(x,y) = \frac{1}{J(x)} I_\Delta (x,y).$

\end{example}

Consider a measured groupoid $(G,\nu)$,
$\nu \in \mathcal{E}$, then, given two functions $\nu$-integrable  $f,g :G\to \mathbb{R}$, we  had defined before an algebra structure on $\mathcal{F}_\nu (G)$ in such way that $(f \,  \underset{\nu}{*}\, g)=h$, if
$$(f \,  \underset{\nu}{*}\, g)(x,y)=
\int \,g(x,s ) \,f (s,y) \nu^{y} (ds) =h(x,y),$$
where $(x,y)\in G$ and $ (s,y)\in G$.

\bigskip

To define the von Neumann algebra associated to $(G,\nu)$, we work with complex valued functions $f:G \to \mathbb{C}$. The product is again given by the  formula

$$(f \,  \underset{\nu}{*}\, g)(x,y)=
\int \,g(x,s ) \,f (s,y) \nu^{y} (ds) .$$

The involution operation $*$ is the rule $f \to \tilde{f}=  f^*,$ where $\tilde{f}(x,y)= \overline{f(y,x)}.$ The functions $f \in \mathcal{F}(G^0)$ are of the form $f(x)=f(x,x)$ are such that $\tilde{f}=f.$

Following Hanh \cite{Hahn2}, we define the I-norm
\[\|f\|_I=  \max \left\{\left\|y\mapsto \int  |f(x,y)| \,\nu^y(dx)\right\|_{\infty},\,\left\|y\mapsto \int  |f(y,x)| \,\nu^y(dx)\right\|_{\infty}\right\},\]
and the algebra $I(G,\nu)=\{f\in L^1(G,\nu):\|f\|_I < \infty\}$ with the product and involution as above. An element $f\in I(G,\nu)$ defines a bounded operator $L_f$ of left convolution multiplication by a fixed $f$ on $L^2(G,\nu)$. This gives the left regular representation of $I(G,\nu)$.

\begin{definition}
Given a measured groupoid $(G,\nu)$, we define the {\bf von Neumann Algebra associated to $(G,\nu)$}, denoted by $W^*(G,\nu)$, as the the von Neumann generated by the left regular representation of $I(G,\nu)$, that is, $W^*(G,\nu)$ is the closure of $\{L_f:f\in I(G,\nu)\}$ in the weak operator topology.

The  multiplicative unity is denoted by $\mathfrak{1}$.

\medskip

In the case $\nu$ is such that such that,  $\int \nu^{y} (ds)=1$, for any $y \in G^0$, we say that the von Neumann algebra is normalized.

\end{definition}

In the setting of  von Neumann Algebras we do not require that $\mathfrak{1}$ is continuous.

\medskip

\begin{definition} We say an element $h\in W^* (G\,,\nu) $ is positive if there exists a $g$ such that $h=g * \tilde{g}$.

This means
$$h(x,y)= (g \,  \underset{\nu}{*}\, \tilde{g})(x,y)=
\int \,g(x,s ) \,\overline{g (y,s)} \nu^{y} (ds) =h(x,y).$$

\end{definition}

Note que $h(x,x) = (g \,  \underset{\nu}{*}\, \tilde{g})(x,x)\geq 0.$

\begin{example} \label{sim} Consider over the set $G^0 =\{1,2..,d\}$ the equivalence relation where all points are related. In this case $G  =\{1,2..,d\} \times \{1,2..,d\}.$ Take $\nu$ as the counting measure. A function $f: G \to \mathbb{R}$ is denoted by $f(i,j),$ where $i\in \{1,2..,d\}, j \in \{1,2..,d\}.$

The convolution product is
$$(f \,  \underset{\nu}{*}\, g)(i,j)= \sum_k g(i,k)\,f(k,j) .$$

In this case the associated von Neumann algebra (the  set of functions $f: G \to \mathbb{C}$)
is identified with the set of matrices, the convolution is the product of matrices and the
identity matrix is the unit $\mathfrak{1}.$ The involution operation is to take the hermitian of a matrix.
\end{example}

\medskip

\begin{example} \label{o1}
For the groupoid $G$ of Example \ref{exi} and the counting measure,   given $f,g :G \to \mathbb{C}$, we have that
$$(f \,  \underset{\nu}{*}\, g)(x,y) =$$
$$\sum_{a\in \{1,2,..,d\}} g(\,(\,x_1,x_2,...\,)\,, \,(\, a,x_2,x_3,..)\,)   \,\,f(\,(\, a,x_2,x_3,..)\,, \,(\,y_1,x_2,...\,)\,) .$$

 We call standard von Neumann algebra on the groupoid $G$ (of Example \ref{exi})
 the associated von Neumann algebra. For this $W^* (G\,,\nu) $ the neutral element $\mathfrak{1}$ (or, more formally $L_\mathfrak{1}$)
is the indicator function of the diagonal (a subset  of $G$). In this case $\mathfrak{1}$
 is measurable but not continuous.

\medskip

\end{example}

\begin{example} \label{J1}
For the probabilistic Haar system $(G,\nu)$ of Example \ref{exi2},   given $f,g :G \to \mathbb{C}$, we get
$$(f \,  \underset{\nu}{*}\, g)(x,y) =$$
$$\sum_{a\in \{1,2,...,d\}} \varphi(  a,x_2,x_3,... ) \,g((\,x_1,x_2,...\,), (\, a,x_2,x_3..))\,\,  f((\, a,x_2,x_3,...), (\,y_1,x_2,...\,)),$$
where $\varphi$ is Holder and such that $\sum_{a\in \{1,2,...,d\}} \varphi(  a,x_1,x_2,... )=1$, for all $x=(x_1,x_2,...)$.

This $\varphi$ is a Jacobian.

The neutral element is described  in example \ref{exi3}.

\end{example}

\begin{example} \label{J0}
In the case $ \nu^y= \delta_{x_0}$ for a fixed $x_0$ independent of $y$, then
$$(f \,  \underset{\nu}{*}\, g)(x,y)=  g(x,x_0)\,f(x_0,y) .$$
\end{example}

\medskip

\begin{proposition} \label{prac}

If $(G,\nu)$ is a measured groupoid, then for $f,g\in I(G,\lambda)$.
 $$(f \,  \underset{\nu}{*}\, g)^{\sim}= \tilde{g} \, \underset{\nu}{*}  \, \tilde{f}. $$

\end{proposition}

 {\bf Proof:}
Remember that for $(x,y)$ in $G$
$$(f \,  \underset{\nu}{*}\, g)(x,y)=
\int \,g(x,s ) \,f (s,y) \nu^{y} (ds) =h(x,y).$$

Then,
$$(f \,  \underset{\nu}{*}\, g)^{\sim} (x,y)=
\int \,\overline{g(y,s ) \,f (s,x)} \nu^{x} (ds).$$

On the other hand
$$(\tilde{g} \,  \underset{\nu}{*}\, \tilde{f})(y,x)=
\int \,\overline{f(s,x)} \,\overline{g (y,s)} \nu^{y} (ds) .$$

As $\nu^{y}= \nu^{x}$ we get that the two expressions are equal.

\qed

\medskip

Then by proposition \ref{prac} we have for the involution $*$ it is valid the property

$$(f \,  \underset{\lambda}{*}\, g)^*= g^* \, \underset{\lambda}{*}  \, f^*. $$

For more details about properties related to this definition we refer the reader to chapter II in \cite{Ren0} and section 5 in \cite{Kas}.

\medskip

We say that
$c:G \to \mathbb{R}$ is a linear cocycle function if $c(x,y) + c(y,z)= c(x,z), $ for all $x,y,z$ which are related. If $c$ is a linear cocycle then $e^\delta$ is a modular function (or, a multiplicative cocycle).

\medskip

\begin{definition} \label{aut}   Consider the von Neumann algebra $W^*(G,\nu)$  associated to $(G,\nu)$.

Given a continuous cocycle function $c:G \to \mathbb{R}$ we define the {\bf group homomorphism} $\alpha:\mathbb{R} \to  \text{Aut} (W^*(G,\nu))$, where for each  $t\in \mathbb{R}$ we have that
$\alpha_t \in $ Aut$(W^*(G,\nu))$ is defined by: for each fixed $t\in  \mathbb{R}$ and $f:G \to \mathbb{R}$ we set
$\alpha_t (f)= e^{t\,i\,c}\, f $.

\end{definition}

{\bf Remark:} Observe that in the above definition that for each fixed $t\in  \mathbb{R}$ and any $f:G^0  \to \mathbb{R}$, we have
	$\alpha_t (f)=f$, since $c(x,x)=0$ for all $x\in G^0$.

\bigskip

We are particularly interested here  in the case where $ G^0=\Omega$ or $G^0 = \hat{\Omega}.$

The value $t$ above is related to temperature and not time. We are later going to consider complex numbers $z$ in place of $t$. Of particular interest is $z= \beta i$ where
$\beta $ is related to the inverse of temperature in Thermodynamic
Formalism (or, Statistical Mechanics).

\medskip

\begin{definition} \label{KMSdef} Consider the von Neumann Algebra $W^*(G,\nu)$ with unity $\mathfrak{1}$ associated to $(G,\nu)$.
A  von Neumann {\bf dynamical state} is a linear functional $w$ (acting on the linear space  $W^*(G,\nu)$) of the form $w:W^*(G,\nu) \to \mathbb{C}$, such that, $w(a)\geq 0$, if $a$ is a positive element of $W^*(G,\nu)$, and $w( \mathfrak{1})=1$.

\end{definition}

\begin{example} \label{co1} Consider over $\Omega=\{1,2,..,d\}^\mathbb{N}$ the equivalence relation
 $R$ of Example \ref{exi}  and the Haar system $(G,\nu)$ associated to the counting measure in each fiber $r^{-1} (x) = \{(a,x_2,x_3,...)\,| a \in \{1,2,...,d\}  \},$ where $x=(x_1,x_2,...)$.

Given a probability $\mu$ over $\Omega$ we can define a   von Neumann dynamical state $\varphi_\mu$  in the following way: for $f:G \to \mathbb{C}$ define
 \begin{equation} \label{ai} \varphi_\mu (f)\,=\,\int f(x,x)\, d\mu(x)= \int \, \,f((\,x_1,x_2,x_3,...\,), (\, x_1,x_2,x_3,...))\,\,  d \mu(x).
 \end{equation}

If $h$ is positive, that is, of the form $h(x,y)=\int \,g(x,s ) \,\overline{g (y,s)} \nu^{y} (ds)$, then
$$ \varphi_\mu (h)\,=  \int (\,\,\int \,\|g(x,s )\|^2  \nu^{x} (ds)\,\,)\,\, d \mu(x)\geq 0.$$

Note that $\varphi_\mu \mathfrak{1}=1$.

Then, $\varphi_\mu$
is indeed a   von Neumann dynamical state.

In this case given $f,g:G \to \mathbb{C}$
$$ \varphi_\mu (f \,  \underset{\nu}{*}\, g)\,=\,$$
$$ \int \,
\sum_{a\in \{1,2,..,d\}} \,f((\,x_1,x_2,...\,), (\, a,x_2,x_3..))\,\,  g((\, a,x_2,x_3..), (\,x_1,x_2,...\,))\, d \mu(x). $$

\medskip

\end{example}

\medskip It seems natural to try to obtain dynamical states from probabilities $M$ on $G^0$ (adapting the reasoning of the above example). Then, given a cocycle $c$ it is also natural to ask: what
we should assume on $M$ in order to get a KMS state for $c$?

\medskip

\begin{example} \label{simon} For the von Neumann algebra of complex matrices of example \ref{sim}
taking $p_1,p_2,..,p_d\geq 0$, such that $p_1+p_2+..+p_d=1$, and $\mu= \sum_{j=1}^d \delta_j,$
we consider $\varphi_\mu$ such that
$$\varphi_\mu (A) = A_{11} p_1 + A_{22} p_2+...+ A_{dd} p_d ,$$
where $A_{ij}$ are the entries of $A$.

Note first that $\varphi_\mu (I)=1$.

If $B= A\, A^*$, then the entries $B_{jj}\geq 0$, for $j=1,2,...,d.$

Therefore, $\varphi_\mu$ is a dynamical state on this von Neumann algebra.

\end{example}

\medskip

\begin{example} \label{co2}  Consider over $\Omega=\{1,2,..,d\}^\mathbb{N}$ the equivalence relation
 $R$ of Example \ref{exi2}  and the associated probability Haar system $\nu$.

Given a probability $\mu$ over $\Omega$ we can define a  von Neumann dynamical state $\varphi_\mu$ in the following way: given $f:G \to \mathbb{C}$ we get $\varphi_\mu(f) = \int f(x,x) \,J(x)\,d \mu(x).$ In this way given $f,g$ we have

$$ \varphi_\mu (f \,  \underset{\nu}{*}\, g)\,=\,$$
$$ \int \,
\sum_{a\in \{1,2,..,d\}} J(  a,x_2,.. )\,g((\,x_1,x_2,...\,), (\, a,x_2,..))\,\,  f((\, a,x_2,..), (\,x_1,x_2...\,))\,J(x) d \mu(x). $$

For the neutral multiplicative element $\mathfrak{1}(x,y) = \frac{1}{J(x)} I_\Delta (x,y)$ we get
$$ \varphi_\mu (\mathfrak{1} )\,=\,\int  \frac{1}{J(x)} I_\Delta (x,x)\, J(x) d \mu(x)=1. $$

\end{example}

\medskip

Consider $G$ a groupoid and a  von Neumann Algebra $W^*(G,\nu)$, where $\nu$ is a transverse function,  with the algebra product
$f  \underset{\nu}{*} g$ and involution $f \to \tilde{f}$.

Given a continuous cocycle $c:G \to \mathbb{R}$ we consider $\alpha:\mathbb{R} \to  \text{Aut} (W^*(G,\nu))$, $t \mapsto  \alpha_t$, the associated  homomorphism according to definition \ref{aut}: for each fixed $t\in  {\bf R}$ and $f:G \to \mathbb{R}$ we set
$\alpha_t (f)= e^{t\,i\,c}\, f $.

\medskip

\begin{definition}
An element $a\in W^*(G,\nu)$ is said to be {\bf analytical} with respect to $\alpha$ if the map $t\in \mathbb{R}\mapsto \alpha_t(a)\in W^*(G,\nu)$ has an analytic continuation to the complex numbers.

More precisely, there is a map $\varphi:\mathbb{C}\to W^*(G,\nu)$, such that, $\varphi(t)=\alpha_t (a)$, for all $t\in\mathbb{R}$, and moreover, for every $z_0\in\mathbb{C}$, there is a sequence $(a_n)_{n\in\mathbb{N}}$ in $W^*(G,\nu)$, such that, $\varphi(z)= \sum_{n=0}^\infty  (z-z_0)^n a_n$ in a neighborhood of $z_0$.

\end{definition}

The analytical elements are dense on the von Neumann algebra (see \cite{Ped}).

\medskip

\begin{definition}
 We say that a  von Neumann  dynamical state  $w$ is a {\bf KMS state for $\beta$  and $c$} if
 $$w \,(\,b\,  \underset{\nu}{*}\,(\alpha_{i\, \beta}(a)\,)\,) = w\,(a\, \underset{\nu}{*}b\,),$$
 for any $b$ and any analytical element $a$.

\end{definition}

It follows from general results (see \cite{Ped}) that it is enough to verify: for any $f,g\in I(G,\nu)$ and $\beta \in \mathbb{R} $ we get
\begin{equation} \label{putz}
w \,(\,g\,  \underset{\nu}{*}\,\alpha_{ \beta\,i}(\,f\,)\,) = w \,(\,g\,  \underset{\nu}{*}\,(\,e^{- \beta\,c}\,f\,)\,) = w\,(f\, \underset{\nu}{*}\,g\,).
\end{equation}

\medskip

Consider the functions
$$u(x,y) = (f * g )(x,y) = \int g(x,s) \,f(s,y) \,\nu^y (ds),$$
and

$$ v(x,y) = (\,g * (e^{- \beta\, c}\,f)\,)\,(x,y) = \int e^{-\beta c(x,s)} f(x,s)\,g(s,y)\, \nu^y (ds).$$
\bigskip

Equation (\ref{putz}) means
\begin{equation} \label{pupu}
w(\, u(x,y)\,) = w (\,v(x,y)\,).
\end{equation}

Note that  equation (\ref{putz}) implies that a KMS von Neumann (or, $C^*$)-dynamical state  $w$ satisfies:

a) for any  $f: G^0 \to \mathbb{C}$ and $g:G \to \mathbb{C}:$
\begin{equation} \label{ll}
w \,(\,g\, \underset{\nu}{*}\,f\,) = w\,(f\,\underset{\nu}{*}\,g\,).
\end{equation}

This follows from the fact that for any
$t\in  \mathbb{R}$ and any $f:G^0  \to \mathbb{R}$, we have that
$\alpha_{t}  (f)=f$.

b) if the function $\mathfrak{1}$ depends just on $x\in G^0$, then, for any $\beta$
$$\alpha_{i\, \beta}( \mathfrak{1} )=\mathfrak{1}.$$

c) $w $ is invariant for the group $\alpha_t$, $t\in \mathbb{R}$.  Indeed,
$$ w( \,\alpha_t (f))=w( \mathfrak{1} \underset{\nu}{*}\,\alpha^t (f)) = w ( f \underset{\nu}{*}\,\mathfrak{1})\,= \, w(f).$$

\bigskip

\begin{example} \label{co21}
For the von Neumann algebra ($C^*$-algebra) of complex matrices of examples \ref{sim} and \ref{simon} consider the dynamical evolution $\sigma_t= e^{i\, t \,H}$, $t \in \mathbb{R}$, where $H$ is a diagonal matrix with entries the real numbers $H_{11}=U_1,H_{22}=U_2,...,H_{dd}=U_d$.
The  KMS state $\rho$ for $\beta$ is
$$\rho(A) = A_{1\,1} \rho_1 + A_{2\,2} \rho_2 +...+ A_{d\,d} \rho_d,$$
where
$ \rho_i= \frac{e^{-\beta U_i}}{\sum_{j=1}^d e^{-\beta U_j}},$ $ i=1,2,...,d$, and $A_{i,j}$, $i,j=1,2,...,d$, are the entries of the matrix $A$ (see \cite{Ren2}).

The probability $\mu$ of example  \ref{co1} corresponds in some sense to the probability $\mu=(\rho_1,\rho_2,..,\rho_d)$ on $\{1,2,...,d\}.$ That is, $\rho=\varphi_\mu.$

This is a clear indication that the $\mu$ associated to the KMS state has in some sense a relation with Gibbs probabilities. This property will appear more explicitly on Theorem \ref{poa1} for the case of the bigger than two equivalence relation.

\end{example}

\medskip

Remember that if $c$ is a cocycle, then $c(x,z)= c(x,y) + c(y,z)$, $\forall x \sim y \sim z$, and, therefore,
$$ \delta(x,y) = e^{\,\beta c (x,y)}=   e^{\,-\beta c (y,x)}$$
is a modular function.

\begin{definition} \label{dede} Given a cocycle $c:G \to \mathbb{R}$ we say that  a probability $M $ over $G^0$ satisfies the $(c,\beta)$-KMS condition for the groupoid $(G,\, \nu)$,  if
for any $h\in I(G,\nu)$, we have
\begin{equation} \label{oip} \int \, \int h(s,x) \nu^x (ds) d M(x)= \int \, \int h(x,s) e^{-\beta \,c (x,s)}\,\nu^x (ds) d M(x),
\end{equation}
where $\beta \in \mathbb{R}$.

In this case we will say that $M$ is a {\bf KMS probability}.

\end{definition}

\medskip

The above means that {\bf $M$ is quasi-invariant for $\nu$ and $\delta(x,s)=e^{-\,\beta c(s,x)}.$}

\medskip
When $\beta=1$ and $c$ is of the form $c(s,x)=V(x)-V(s)$ the above condition means
\begin{equation} \label{oipa} \int \, \int h(s,x) \nu^x (ds) e^{V(x)} d M(x)= \int \, \int h(x,s) e^{V(x)} \,\nu^x (ds) d M(x).
\end{equation}

\medskip

\begin{proposition} \label{war5} (J. Renault - Proposition II.5.4 in \cite{Ren0}) Suppose that the
state $w$ is such that for a certain  probability $\mu$  on $G^0$ we have that for any $h\in I(G,\nu)$ we get $w(h) = \int h(x,x) d \mu(x)$. Then, to say that $\mu$ satisfies the $(c,\beta)$-KMS condition for $(G,\nu)$ according to Definition \ref{dede} is equivalent to say  that $w$ is KMS for $(G,\nu)$, $c$ and $\beta$, according to equation \eqref{putz}.
\end{proposition}

{\bf Proof:} Note that  for any $f, g$
$$(\,f\,  \underset{\nu}{*}\,\,g\,)(x,y)= \int g(x,s)\, f(s,y) d \nu^x(s)$$
and
$$(\,g\,  \underset{\nu}{*}\,(\,e^{- \beta\, c}\,f\,)\,)(x,y)= \int f(x,s)\, g(s, y) \,e^{- \beta\, c(x,s)} \,d \nu^x(s). $$

We have to show that $\int u(x,x) d \mu(x) = \int v(x,x) d \mu (x)$ (see equation (\ref{pupu})).

Then, if the $(c,\beta)$-KMS condition for $M$ is true, we take $h(s,x) =  g(x,s)\, f(s,x)$ and we got equation \eqref{pupu} for such $w$.

By the other hand if  \eqref{pupu} is true for such $w$ and any $f,g$, then take $f(s,x) =h(s,x)$ and $g(s,x)=1$.

\qed

\medskip

\begin{example} \label{rr} In the case for each $y$ we have that $\nu^y$ is the counting measure we get that to say that
a probability $M $ over $\hat{\Omega}$ satisfies the $(c,\beta)$-KMS condition means: for any $h:G\to \mathbb{C}$
\begin{equation} \label{ulu} \sum_{y \sim x} \int \, h(x,y) e^{-\beta \,c (x,y)} d M(x)= \sum_{x \sim y}  \, \int h(x,y) \, d M(y).
\end{equation}

In the notation of \cite{Ren1} we can  write the above  in an equivalent way as
$$ \int \, h \,e^{-\beta \,c } d (s^* (M)) =   \, \int h\, d (r^* (M)).
$$

Note that in \cite{Ren1} it is considered $r(x,y)=x$ and $s(x,y)=y.$
\medskip


Suppose  $c(x,y) = \varphi(x) - \varphi(y).$
Then, taking $h(x,y)=k(x,y)\, e^{\beta\,\varphi(x)}$ we get an equivalent expression for (\ref{ulu}): for any $k(x,y)$
\begin{equation} \label{uzu} \sum_{y \sim x} \int \, k(x,y) e^{\beta \,\varphi (y)} d M(x)= \sum_{x \sim y}  \, \int k(x,y) \, e^{\beta \,\varphi (x)}d M(y).
\end{equation}



\end{example}

 For a  Holder continuous potential $A:  \{1,2..,d\}^\mathbb{N} \to \mathbb{R}$ the Ruelle operator $\mathcal{L}_A$  acts on continuous functions $v: \{1,2..,d\}^\mathbb{N} \to \mathbb{R}$ by means of $\mathcal{L}_A(v)=w$, if
 $$\mathcal{L}_A(v)(x_1,x_2,x_3,...) = \sum_{a=1}^d e^{ A(a,x_1,x_2,x_3,...)} v(a,x_1,x_2,x_3,...) = w(x).$$

For a  Holder continuous potential $A:  \{1,2..,d\}^\mathbb{N} \to \mathbb{R}$ there exist a continuous positive eigenfunction $f$, such that, $\mathcal{L}_A (f)=  \lambda \, f$, where $\lambda$ is positive and also the spectral radius of $\mathcal{L}_A$ (see \cite{PP}).

The dual $\mathcal{L}_A^*$  of  $\mathcal{L}_A$ acts on probabilities by Riesz Theorem (see \cite{PP}). We say that the probability $m$ on $\{1,2..,d\}^\mathbb{N}$ is {\bf Gibbs for the potential $A$}, if $\mathcal{L}_A^*(m)= \lambda \,m$ (same $\lambda$ as above). In this case we say that $m$ is an eigenprobability for $A$.

Gibbs probabilities  for Holder potentials $A$ are also {\bf DLR probabilities} on $\{1,2,...d\}^\mathbb{N}$ (see \cite{CL1}).

Gibbs probabilities for Holder potentials $A$ can be also obtained via Thermodynamic Limit from boundary conditions (see \cite{CL1}).

We say that the {\bf potential $A$ is normalized } if $\mathcal{L}_A(1)=1$.
In this case a probability $\mu$ is Gibbs (equilibrium) for the normalized potential $A$ if it is a fixed point for the dual
of the Ruelle operator, that is, $\mathcal{L}_A^* (\mu)=\mu.$

Suppose $\Omega =\{-1,1\}^\mathbb{N}$ and $A: \Omega \to \mathbb{R}$ is of the form
$$ A(x_0,x_1,x_2,...) = x_0 a_0 + x_1 a_1 + x_2 a_2 +x_3 a_3 +...+ x_n a_n+...$$
where $\sum a_n$ is absolutely convergent.

In \cite{CDLS} the explicit
expression of the  eigenfunction for $\mathcal{L}_A$ and the eigenprobability for the dual  $\mathcal{L}_A^*$ of
the Ruelle operator $\mathcal{L}_A$  is presented. The eigenprobability   is not invariant for the shift.

\medskip

\begin{example} \label{win593}

For the Haar system of examples \ref{win590} and  \ref{win591} where $k$ is fixed consider a normalized potential Holder $A:\{1,2..,d\}^\mathbb{N} \to \mathbb{R}$. Denote by $\mu$ the  equilibrium probability associated to such $A$.

Consider $$ \delta(x,y)= \frac{e^{A(y) + A(\sigma(y))+...+ A (\sigma^{k-1}(y))}}{e^{A(x) + A(\sigma(x)+...+ A (\sigma^{k-1}(x))}}.$$

We claim that $\mu$ satisfies the
$(c,\beta)$-KMS condition  (\ref{oipa}) for such $\delta$ when $\beta=1$.

For each cylinder set $\overline{a_1,a_2,..,a_k}$ the transformation $\sigma^k : \overline{a_1,a_2,..,a_k} \to \{1,2..,d\}^\mathbb{N}$ is a bijection.  The pull back by $\sigma^k$  of the probability $\mu$ with respect to $\mu$ has Radon-Nykodin derivative
$$\phi_{ \overline{a_1,a_2,..,a_k} } (x)= e^{A(x) + A(\sigma(x)+...+ A (\sigma^{k-1}(x))}.$$

Denote $\varphi_{ \overline{a_1,a_2,..,a_k} }\,:\, \{1,2..,d\}^\mathbb{N} \to \overline{a_1,a_2,..,a_k} $ the inverse of $\sigma^k $ (restricted to $\{1,2..,d\}^\mathbb{N}$).

Consider the cylinders $K=\overline{a_1,a_2,..,a_k}$, $L=\overline{b_1,b_2,..,b_k}$.

Note that it follows from the use of the change of coordinates $y \to x=   \varphi_{ \overline{b_1,b_2,..,b_k} } \circ (\,\varphi_{ \overline{a_1,a_2,..,a_k} }\,)^{-1} (y)$ that
$$ \int_K  e^{A(y) + A(\sigma(y))+...+ A (\sigma^{k-1}(y))}  d \mu(y)= \int_L e^{A(x) + A(\sigma(x)+...+ A (\sigma^{k-1}(x))}\, d\mu(x)$$

For each class the number of elements $s$ on $K$ or $L$ is the same.

This means that
$$  \sum_{s\in L}\, \int_K    e^{A(y) + A(\sigma(y))+...+ A (\sigma^{k-1}(y))}  d M(y)=$$
\begin{equation} \label{kwe59}
 \, \sum_{s\in K}\int_L  \,e^{A(x) + A(\sigma(x))+...+ A (\sigma^{k-1}(x))}  d M(x) .
\end{equation}

\medskip

{\bf Remark:} Given $y\in X$, consider  the function $f(x,s),$ where
$f: [y] \times [y] \to \mathbb{C}.$

For each pair $(i,j) \in [y] \times [y]$, denote $z_{i,j} = f(i,j).$

Then, $f: [y] \times [y] \to \mathbb{C}$ can be written as
$$ \sum_{i,j \in \, [y]} z_{i,j} I_{i}\, I_{j}.$$

Then, any  function $f(x,s)$, $f: [y] \times [y] \to \mathbb{C},$ is a linear combination of functions which are the product of two functions: one depending just on $x$ and the other just on $s$.

\medskip

Then, expression  (\ref{kwe59})
means that for such $f$ we have
$$  \sum_s\, \int f(s,y) e^{A(y) + A(\sigma(y))+...+ A (\sigma^{k-1}(y))}  d M(y)=$$
\begin{equation} \label{kwe37}
 \sum_s\, \int  f(x,s) e^{A(x) + A(\sigma(x))+...+ A (\sigma^{k-1}(x))}  \, d M(x),
\end{equation}

The $(c,\beta)$-KMS condition  (\ref{oipa}) for the probability  $M$  and for any continuous  function $f$ means
$$  \sum_s\, \int f(s,y) e^{\beta (A(y) + A(\sigma(y))+...+ A (\sigma^{k-1}(y))}  d M(y)=$$
\begin{equation} \label{kwe58}
 \sum_s\, \int  f(x,s) e^{\beta (A(x) + A(\sigma(x))+...+ A (\sigma^{k-1}(x)))}  \, d M(x),
\end{equation}

Expression (\ref{kwe58})  follows  from (\ref{kwe37}) and the above remark. Therefore, such $M$ satisfies the KMS condition for such $\delta$.

\end{example}

\medskip

 In example \ref{o1} consider  $\Omega=\{1,2\}^\mathbb{N}$ and take $\nu^y$ the counting
measure on the class of $y$.
Consider the von Neumann algebra associated to this measured groupoid $(G,\nu)$ where $G$ is given by the bigger than two relation.

In this case $\mathfrak{1}(x,y) = I_\Delta (x,y)$.

Consider  $c(x,y) = \varphi(x) - \varphi(y)$, where $\varphi$ is Holder.
We do    not assume that  $\varphi$ is normalized.

 A natural question is: the eigenprobability  $\mu$ for such potential
$\varphi$ is such that $f \to \varphi_\mu(f)\,=\,\int f(x,x) \,\,d \mu(x)$ defines the associated KMS state?
For each modular function $c$?

The purpose of the next results is to analyze this question when $c(x,y) = \varphi(x) - \varphi(y).$

\medskip

Consider the equivalence relation on $\Omega =\{1,2...,d\}^\mathbb{N}$ which is
$$ x=(x_1,x_2,x_3,..) \sim y=(y_1,y_2,y_3,...)  \,\,,\, \text{if an only if } \,,\, x_j=y_j\,\,\,\text{ for all}\,\,\, j\geq 2.$$

In this case the class $[x]$  of  $ x=(x_1,x_2,x_3,..) $ is
$$[x] = \,\{ \,(1,x_2,x_3,..),\, (2,x_2,x_3,..),..., (d,x_2,x_3,..)\,\}.$$

The associated groupoid by $G\subset \Omega \times \Omega$, is
$$  G= \{(x,y)\,|\, x \sim y\}.$$

$G$ is a closed set on the compact set $\Omega \times \Omega$.
We fix the measured groupoid $(G,\nu)$ where $\nu^x$ is the counting measure. The results we will get are the same if we
take the Haar system as the one where each point $y$ on the class of $x$ has mass $1/d$.

In this case equation (\ref{oip}) means
$$
\sum_j \, \int f(\,(j,x_2,x_3,..,x_n,..) \,,\,(x_1,x_2,x_3,..,x_n,..)  ) d M(x)=$$
\begin{equation} \label{oipmi} \sum_j \, \int f(\,(x_1,x_2,x_3,..) \,,\,(j,x_2,x_3,..)  )
e^{ - c(j,x_2,x_3,..) \,,\,(x_1,x_2,x_3,..)  ) } d M(x).
\end{equation}

The first question: given a cocycle $c$ does there exist $M$ as above?

\medskip

Suppose
$c(x,y) = \varphi(y) - \varphi(x)$.

In this case equation (\ref{oipmi}) means
$$
\sum_j \, \int f(\,(j,x_2,x_3,..,x_n,..) \,,\,(x_1,x_2,x_3,..,x_n,..)  ) d M(x)=$$
\begin{equation} \label{oipal} \sum_j \, \int f(\,(x_1,x_2,x_3,..) \,,\,(j,x_2,x_3,..)  )
e^{ - \varphi (j,x_2,x_3,..) \,+\,\varphi (x_1,x_2,x_3,..)   } d M(x).
\end{equation}

 Among other things we will show later that if we assume that  $\varphi$ depends just on the first coordinate then we can take $M$ as the independent probability (that is, such independent $M$ satisfies the KMS condition (\ref{oipal})).

\medskip

In section 3.4 in \cite{Ren2} and in  \cite{KumRen} the authors present a result concerning quasi-invariant probabilities and Gibbs probabilities on $\{1,2,..d\}^\mathbb{N}$ which has a different nature when compared to the next one. The groupoid is different from the one we will consider (there elements are of the form $(x,n,y)$, $n \in \mathbb{Z}).$
In \cite{Ren2} and \cite{KumRen}  for just one value of $\beta$ you get the existence of the quasi invariant probability. Moreover, the KMS state is unique (here this will not happen as we will show on Theorem \ref{KMS1})

 In \cite{Ru1}, \cite{HD}, \cite{Meye} and \cite{LM2}  the authors present  results which have some similarities with the next theorem. They consider  Gibbs  (quasi-invariant) probabilities in the case of the symbolic space $\{1,2,..d\}^\mathbb{Z}$ and not $\{1,2,..d\}^\mathbb{N}$ like here. In all these papers the quasi-invariant probability is unique and invariant for the shift. In  \cite{Bis} the authors consider DLR probabilities for interactions in $\{1,2,..d\}^\mathbb{Z}$. The equivalence relation (the homoclinic relation of Example \ref{homore}) in all these cases is quite different from the one we will consider.

\medskip

\begin{theorem} \label{poa1}
Consider the Haar system with the counting measure $\nu$ for the bigger than two relation on $\{1,2,..d\}^\mathbb{N}$.
Suppose  that $\varphi$ depends just on the first $k$  coordinates, that is,
and
$$\varphi(x_1,x_2,...,x_k, x_{k+1} ,x_{k+2},..)= \varphi(x_1,x_2,...,x_k).$$
Then, the eigenprobability $\mu$ (a DLR probability) for the potential $-\varphi$  (that is,
$\mathcal{L}_{-\varphi}^* (\mu)=\lambda \mu,$ for some positive $\lambda$)
satisfies the KMS condition (is quasi-invariant) for the associated modular function  $c(x,y) = \varphi(y)- \varphi(x)$.

The same result is true, of course, for $\beta c$, where $\beta>0.$

\end{theorem}

{\bf Proof:} We are going to show that the Gibbs probability $\mu$ for the potential $-\varphi$  satisfies the KMS condition.

We point out that in general the eigenvalue $\lambda \neq 1$.

We have to show that (\ref{oipal}) is true when $M=\mu$. That is, $\mu$ is a KMS probability for the Haar system and the modular function.

Denote for any finite string $a_1,a_2,...a_n $ and any $n$

$$ p_{ a_1,a_2,...a_n  } =\frac{ e^{-[\varphi ( a_1,a_2,...a_n, \,1^\infty) +  \varphi( a_2,...a_n, \,1^\infty ) +...  +   \varphi( a_n, \,1^\infty)] } }{ \sum_{b_1,b_2,...b_n} e^{-\,[\varphi ( b_1,b_2,...b_n, \,1^\infty) +  \varphi( b_2,...b_n, \,1^\infty ) +...  +   \varphi( b_n, \,1^\infty)] } }     .$$

Note that for $n>k$ we have that
$$e^{-[\,\varphi ( a_1,a_2,...a_n, \,1^\infty) +  \varphi( a_2,...a_n, \,1^\infty ) +...  +
\varphi( a_n, \,1^\infty) \,]}=$$
$$e^{-[\,\varphi ( a_1,a_2,...a_k) +  \varphi( a_2,...a_{k+1} ) +...  +   \varphi( a_n,\underbrace{ 1,1,..,1}_{k-1} )+ (n-k) \varphi( \underbrace{ 1,1,..,1}_{k-1}  )] }.$$

Therefore,
$$ \sum_{b_1,b_2,...b_n} e^{-\,[\,\varphi ( b_1,b_2,...b_n, \,1^\infty) +  \varphi( b_2,...b_n, \,1^\infty ) +...  +   \varphi( b_n, \,1^\infty) ]}=$$
$$ e^{ - (n-k) \varphi( \underbrace{ 1,1,..,1}_{k-1}  )}\,\sum_{b_1,b_2,...b_n} e^{-[\,\varphi ( b_1,b_2,...b_k) +  \varphi( b_2,...b_{k+1} ) +...  +   \varphi( b_n,\underbrace{ 1,1,..,1}_{k-1} ) \,]}.$$

Consider the probability $\mu_n$, such that,

$$\mu_n = \sum_{a_1,a_2,...a_n} \delta_{(a_1,a_2,...a_n, \,1^\infty)} \,  p_{ a_1,a_2,...a_n  }=$$
$$   \sum_{a_1,...a_n} \delta_{(a_1,...a_n, \,1^\infty)} \, \frac{ e^{-[\,\varphi ( a_1,...a_k) +  \varphi( a_2,...a_{k+1} ) +...  +   \varphi( a_n,\underbrace{ 1,..,1}_{k-1} ) ]}}{  \sum_{b_1,...b_n} e^{-[\,\varphi ( b_1,...b_k) +  \varphi( b_2,...b_{k+1} ) +...  +   \varphi( b_n,\underbrace{ 1,..,1}_{k-1} ) ]}}.$$

and $\mu$ such that $\mu = \lim_{n \to \infty} \mu_n.$

Note that
$$ p_{a_1,a_2,..,a_n} =  \frac{ e^{-[\,\varphi ( a_1,...a_k) +  \varphi( a_2,...a_{k+1} ) +...  +   \varphi( a_n,\underbrace{ 1,..,1}_{k-1} ) \,]}}{  \sum_{b_1,...b_n} e^{-[\, \varphi ( b_1,...b_k) +  \varphi( b_2,...b_{k+1} ) +...  +   \varphi( b_n,\underbrace{ 1,..,1}_{k-1} )\,] }}.$$

If $\varphi$ is Holder  it is known that the above probability $\mu$ is the
eigenprobability for the dual of the Ruelle operator $\mathcal{L}_{-\varphi}$ (a DLR probability).
That is, there exists $\lambda>0$ such that $ \mathcal{L}_{-\varphi}^* (\mu)= \lambda \mu$. This follows from the Thermodynamic Limit with boundary condition property as presented in \cite{CL1}.

We claim  that the above probability $\mu$ satisfies the KMS condition.


Indeed,
note that
$$
\sum_j \,  \int f(\,(j,x_2,x_3,..,x_n,..) \,,\,(x_1,x_2,x_3,..,x_n,..)  ) d \mu(x)=$$
$$\lim_{n \to \infty} \sum_j \, \sum_{a_1,a_2,...a_n} \, f(\,(j,a_2,a_3,..,a_n,\, 1^\infty ) \,,\,(a_1,a_2,a_3,..,a_n,\, 1^\infty) \, )   p_{a_1,a_2,...,a_n} = $$


 \begin{equation} \label{lele}\lim_{n \to \infty} \sum_j \, \sum_{a_1}\sum_{a_2,...a_n}   f((j,a_2,a_3,..,a_n,\, 1^\infty ) (a_1,a_2,a_3,..,a_n, 1^\infty) p_{a_1,a_2,...,a_n}     .
 \end{equation}

On the other hand

$$ \sum_j \, \int f(\,(x_1,x_2,x_3,..) \,,\,(j,x_2,x_3,..)  )
e^{ - \varphi (j,x_2,x_3,...) \,+\,\varphi (x_1,x_2,x_3,..)   } d \mu(x)=$$
$$\lim_{n \to \infty}\sum_j \, \sum_{a_1,a_2,...,a_n} \, f(\,(a_1,..,a_n,\, 1^\infty ) \,,\,(j,a_2,..,a_n,\, 1^\infty) \, )  e^{ - \varphi (j,a_2,...,a_k) \,+\,\varphi (a_1,a_2,...a_k)   } p_{a_1,a_2,...,a_n}  =$$
$$\lim_{n \to \infty}\sum_j \,\sum_{a_1} \sum_{a_2,...a_n} \, f(\,(a_1,..,a_n,\, 1^\infty ) \,,\,(j,a_2,..,a_n,\, 1^\infty) \, ) e^{ - \varphi (j,a_2,...,a_k) \,+\,\varphi (a_1,a_2,...a_k)   } $$
 $$\frac{ e^{-\,[\,\varphi ( a_1,a_2,...a_k) +  \varphi( a_2,...a_{k+1} ) +...  +   \varphi( a_n,\underbrace{ 1,..,1}_{k-1} ) ]}}{  \sum_{b_1,...b_n} e^{-[\,\varphi ( b_1,...b_k) +  \varphi( b_2,...b_{k+1} ) +...  +   \varphi( b_n,\underbrace{ 1,..,1}_{k-1} ) ]}}=$$
 $$\lim_{n \to \infty}\sum_j \,\sum_{a_1} \sum_{a_2,...a_n} \, f(\,(a_1,..,a_n,\, 1^\infty ) \,,\,(j,a_2,..,a_n,\, 1^\infty) \, ) $$
 $$\frac{ e^{-\,[\,\varphi ( j,a_2,,...a_k) +  \varphi( a_2,...a_{k+1} ) +...  +   \varphi( a_n,\underbrace{ 1,..,1}_{k-1} ) ]}}{  \sum_{b_1,...b_n} e^{-[\,\varphi ( b_1,...b_k) +  \varphi( b_2,...b_{k+1} ) +...  +   \varphi( b_n,\underbrace{ 1,..,1}_{k-1} ) ]}}.=$$
 $$\lim_{n \to \infty}\sum_j \,\sum_{a_1} \sum_{a_2,...a_n} \, f(\,(a_1,..,a_n,\, 1^\infty ) \,,\,(j,a_2,..,a_n,\, 1^\infty) \, )p_{j,a_2,...,a_n}. $$

On this last equation if we exchange coordinates $j$ and $a_1$ we
get expression (\ref{lele}).

Then, such $\mu$ satisfies the KMS condition.

\bigskip

\qed

The above theorem can be extended to the case the potential $\varphi$ is Holder.
We refer the reader to \cite{LO} for more general results.

\bigskip

We will show now that under the above setting the KMS probability is not unique.

\medskip

\begin{proposition} \label{KMS1} \,Suppose  $\mu $ satisfies the  KMS
condition for the measured groupoid $(G,\nu)$ where
$c(x,y) = \varphi(y) - \varphi(x)$. Suppose $\varphi$ is normalized for the Ruelle operator,
where $\varphi: G^0 =\Omega \to \mathbb{R}$.  Consider $v(x_1,x_2,x_3,..)$ which does not depend of the first coordinate.
Then, $v(x) d\mu(x)$ also satisfies the  KMS
condition for the measured groupoid $(G,\nu)$.

\end{proposition}

{\bf Proof:}
Suppose  $\mu $ satisfies the $(c,\beta)$-KMS
condition for the measured groupoid $(G,\nu)$.
This means: for any $g\in I(G,\nu)$

 $$\int \,
\sum_{a\in \{1,2,..,d\}} \,g((\, a,y_2,y_3..)\,,\,(\,y_1,y_2,...\,) )\,e^{ \beta \varphi(\, a,y_2,y_3..) }\,\, d \mu(y)=$$
\begin{equation} \label{uku3} \int \,
\sum_{a\in \{1,2,..,d\}} \,g((\,x_1,x_2,...\,), (\, a,x_2,x_3..))\,e^{ \beta \varphi(\, a,x_2,x_3..) }\,\, d \mu(x).
\end{equation}

Take
$$h(x_1,x_2,x_3,..), (y_1,y_2,y_3,..) ) = $$
$$k(\,(x_1,x_2,x_3,..), (y_1,y_2,y_3,..) )\, v(x_1,x_2,x_3,..)\,)=$$
$$k(\,(x_1,x_2,x_3,..), (y_1,y_2,y_3,..) )\, v(x_2,x_3,..)\,).$$

From the hypothesis about $\mu$ we get that

 $$\int \,
\sum_{a\in \{1,2,..,d\}} \,h((\, a,y_2,y_3..)\,,\,(\,y_1,y_2,...\,) )\,e^{ \beta \varphi(\, a,y_2,y_3..) }\,\, d \mu(y)=$$
$$ \int \,
\sum_{a\in \{1,2,..,d\}} \,h((\,x_1,x_2,x_3,...\,), (\, a,x_2,x_3..))\,e^{ \beta \varphi(\, a,x_2,x_3..) }\,\, d \mu(x).
$$

This means, for any continuous $k$ the equality
$$\int \,
\sum_{a\in \{1,2,..,d\}} \,k((\, a,y_2,y_3..)\,,\,(\,y_1,y_2,...\,) )\,e^{ \beta \varphi(\, a,y_2,y_3..) }\,\,v(y_2,y_3,...) d \mu(y)=$$
$$ \int \,
\sum_{a\in \{1,2,..,d\}} \,k((\,x_1,x_2,x_3,...\,), (\, a,x_2,x_3..))\,e^{ \beta \varphi(\, a,x_2,x_3..) }\,\,v(x_2,x_3,..) d \mu(x).
$$

Therefore,
$v(x) d\mu(x)$ also satisfies the $(c,\beta)$-KMS
condition for the measured groupoid $(G,\nu)$.

\qed

\medskip
It follows from the above result that the probability that satisfies the KMS
condition for $c$ and  the measured groupoid $(G,\nu)$ is not always unique.

A probability $\rho$ satisfies the Bowen condition for the potential $-\varphi$ if
there exists constants $c_{1},c_{2} > 0$, and P, such that,
for every
$$x = (x_{1},...,x_m,...) \in \Omega = \{1,2,...d\}^\mathbb{N},$$
and all $m \geq 0$,
\begin{equation} \label{defbowen}
c_{1} \leq
\frac{\rho \{y : y_{i} = x_{i}, \quad \forall i = 1, . . . , m \}  }
{\exp \left( -Pm - \sum_{k=1}^{m} \varphi(\sigma^{k}(x) \right) }
\leq c_{2}.
\end{equation}

Suppose $\varphi$ is Holder, then,
if $\rho$ is the equilibrium probability (or, if $\rho$ is the eigenprobability for the
dual of Ruelle operator $\mathcal{L}_{-\varphi}$) one can show that
it satisfies the Bowen condition for $-\varphi$.

 In the case $v$ is continuous and  does not depend on the first coordinate then
$v(x) d\mu(x)$ also satisfies the Bowen condition for $\varphi$. The same is true for the
probability $\hat{\rho}$ of example \ref{oii} on the case $-\varphi= \log J$.

There is an analogous definition of the Bowen condition on the space $\{1,2,...d\}^\mathbb{Z}$ but
it is a much more strong hypothesis on this case (see section 5 in \cite{LM2}).

\bigskip

\begin{example} \label{oii} We will show an example where the probability $\mu$ of  theorem \ref{poa1} (the eigenprobability for the potential $-\varphi$) is such that if $f$ is a function that depends just on the first coordinate, then,
$f \, \mu$ does not necessarily satisfies the KMS condition.

Suppose $\varphi = -\log J $, where $J(x_1,x_2,x_3,..) =  J(x_1,x_2)>0$, and
$\sum_i P_{i,j}=1$, for all $i$. In other words the matrix $P$, with entries $P_{i,j}$,
$i,j \in \{1,2..,d\}$,  is a column stochastic matrix.
The Ruelle operator for $-\varphi$ is the Ruelle operator for $ \log J.$ The potential $\log J$ is normalized for
the Ruelle operator.

We point out that in Stochastic Process it is usual to consider line stochastic matrices which is different
from our setting.

There exists a unique right eigenvalue probability vector $\pi$  for $P$ (acting on vectors on the right). The Markov chain determined by
the matrix $P$ and the initial vector of probability $\pi= (\pi_1,\pi_2,..,\pi_d)$ determines an stationary process, that is,
a probability $\rho$ on the Bernoulli space $\{1,2,...,d\}^\mathbb{N},$  which is invariant for the shift
acting on $\{1,2,...,d\}^\mathbb{N}$.

For example, we have that $\rho(\overline{21})= P_{21} \pi_1$.

We point out that such $\rho$  is  the eigenprobability
for the $\mathcal{L}^*_{\log J}$ (associated to the eigenvalue $1$).
Therefore, $\rho$ satisfies the KMS condition from the above results.

The Markov Process determined by
the matrix $P$ and the initial vector of probability $\pi= (1/d,1/d,...,1/d)$ defines
a probability $\hat{\rho}$ on the Bernoulli space $\{1,2,...,d\}^\mathbb{N},$  which is not invariant for the shift
acting on $\{1,2,...,d\}^\mathbb{N}$.

 In this case, for example, $\hat{\rho}(\overline{21})= P_{21} \,1/d$.

 Note that the probability $\rho$ satisfies  $\rho= u \,\hat{\rho}$ where $u$ depends just on the first coordinate.

Note that unless $P$ is double stochastic is not true that for any $j_0$ we have that $\sum_{k} \,
P_{j_0,k}\,=1$.

Assume that there exists $j_0$ such that  $\sum_{k} \,
P_{j_0,k}\,\neq 1$.

We will check that, in this case  $\hat{\rho}$ does not satisfies the KMS condition for  the function
$f(x,y) =I_{X_1=i_0} (x) \,I_{X_1=j_0} (y).  $

Indeed, equation (\ref{oipal}) means

$$
\sum_j \, \int f(\,(j,x_2,x_3,..,x_n,..) \,,\,(x_1,x_2,x_3,..,x_n,..)  ) d \hat{\rho}(x)=$$
$$\sum_j \int I_{X_1=i_0} (j,x_2,x_3,..) \,I_{X_1=j_0} (x_1,x_2,x_3...) d \hat{\rho} (x) = $$
$$ \int I_{X_1=i_0} (i_0,x_2,x_3,..) \,I_{X_1=j_0} (x_1,x_2,x_3...) d \hat{\rho}(x) = $$
$$  \int  \,I_{X_1=j_0} (x_1,x_2,x_3...) d \hat{\rho} (x)   = \hat{\rho} (\overline{j_0})=\,1/d =$$
$$\sum_j \, \int f(\,(x_1,x_2,x_3,..) \,,\,(j,x_2,x_3,..)  )
e^{ - \varphi (j,x_2,x_3,..) \,+\,\varphi (x_1,x_2,x_3,..)   } d \hat{\rho}(x)=
$$
$$\sum_j \,\int  I_{X_1=i_0} (x_1,x_2,x_3,..) \,I_{X_1=j_0} (j,x_2,x_3...)
e^{ - \varphi (j,x_2,x_3,..) \,+\,\varphi (x_1,x_2,x_3,..)   } d \hat{\rho}(x)=
$$
$$\int  I_{X_1=i_0} (x_1,x_2,x_3,..) \,I_{X_1=j_0} (j_0,x_2,x_3...)
e^{ - \varphi (j_0,x_2,x_3,..) \,+\,\varphi (x_1,x_2,x_3,..)   } d \hat{\rho}(x)=
$$
$$\int  I_{X_1=i_0} (x_1,x_2,x_3,..) \,
e^{ - \varphi (j_0,x_2,x_3,..) \,+\,\varphi (x_1,x_2,x_3,..)   } d \hat{\rho}(x)=
$$
$$\int_{X_1=i_0} \,
e^{ - \varphi (j_0,x_2,x_3,..) \,+\,\varphi (i_0,x_2,x_3,..)   } d \hat{\rho}(x)=
$$
$$\sum_{k} \int_{X_1=i_0\,,\, X_2=k} \,
e^{ - \varphi (j_0,x_2,x_3,..) \,+\,\varphi (i_0,x_2,x_3,..)   } d \hat{\rho}(x)=
$$
$$\sum_{k} \int_{X_1=i_0\,,\, X_2=k} \,
P_{j_0,k}\, \, P_{i_0,k}^{-1} d \hat{\rho}(x)=
$$
$$\sum_{k}  \,
P_{j_0,k}\, \, P_{i_0,k}^{-1} \,  P_{ i_0,k} 1/d=
$$
$$ \sum_{k} \,
P_{j_0,k}\, \,1/d \, \neq 1/d= \hat{\rho}(\overline{j_0}).
$$

Therefore, $\hat{\rho}$ does not satisfies the KMS condition.

\end{example}

\begin{example}

Consider  $\Omega=\{1,2\}^\mathbb{N}$, a  Jacobian $J$ and take $\nu^y$ the
probability  on each class  $y$ given by $\sum_a J(a,y_2,y_3,..) \delta_{(a,y_2,y_3,..)  }.$

Note first that  $\varphi= \log J$ is a  normalized potential. Does the equilibrium probability for $\log J$ satisfies the KMS condition? We will show that this in not always true.

The  question means: is it true that for any function $k$ is valid
$$ \int \,
\sum_{a\in \{1,2\}} k( \,( a,y_2,.. ), \,(y_1,y_2,..)\,)\,e^{ \varphi  (a,y_2,...)}\, d \mu(y)= $$
$$ \int \,
\sum_{a\in \{1,2\}} k( \,( a,x_2,.. ), \,(x_1,x_2,..)\,)\,e^{ \varphi  (a,x_2,...)}\, d \mu(x)= $$
\begin{equation} \label{Kre}
\int \,
\sum_{a\in \{1,2\}} k( \,( x_1,x_2,.. ), \,(a,x_2,..)\,)\,e^{ \varphi (a,x_2,...)}\, d \mu(x)?
\end{equation}

Consider the example:
take  $c(x,y) = \varphi(x) - \varphi(y)$, for $\varphi: \{1,2\}^\mathbb{N} \to \mathbb{R}$,  such that,
$$  \varphi (a,\,.\,,\,.\,,...) = \log p$$
where $p=  p_a$,\,\,for $\,a\in\{1,2\}$,
and $p_1+p_2=1,$ $p_1,p_2>0.$

The Gibbs probability $\mu$ for such $\varphi$ is the independent probability associated to $p_1,p_2$.

Given such probability $\mu$ over $\Omega$ we can define a dynamical state $\varphi_\mu$ in the following way: given $f:G \to \mathbb{R}$ we get $\varphi_\mu(f) = \int f(x,x) \,\,d \mu(x).$

Take $\beta=1$. We will show that $\varphi_\mu$ is not KMS for $c$.

The equation  (\ref{Kre}) for such $\mu$ means for  any $k(x,y)$
$$ \int \,
\sum_{a\in \{1,2\}} p_a\, k( \,( a,y_2,.. ), \,(y_1,y_2,..)\,)\,\, d \mu(y)= $$
$$ \int \,
\sum_{a\in \{1,2\}} k( \,( a,y_2,.. ), \,(y_1,y_2,..)\,)\,p (a,y_2,...)\, d \mu(y)= $$
$$ \int \,
\sum_{b\in \{1,2\}} k( \,( x_1,x_2,.. ), \,(b,x_2,..)\,)\,p (b,x_2,...)\, d \mu(x)= $$
$$ \int \,
\sum_{b\in \{1,2\}} p_b\, k( \,( x_1,x_2,.. ), \,(b,x_2,..)\,)\,\, d \mu(x).$$

It is not true that $\mu$ is Gibbs for the potential $\log p$.

Indeed, given $k$ consider the function
$$g(y_1,y_2,y_3,y_4,...)= k( \,(y_1,y_3,...)\,, \,( y_2,y_3,...)\,).$$

Note that
$$ \mathcal{L}_{\log p} (g)(y_1,y_2,y_3,..)= \sum_{a\in \{1,2\}} p(a,y_1,y_2,y_3,...) \, g(a,y_1,y_2,...)$$
$$\sum_{a\in \{1,2\}} p_a\, k( \,( a,y_2,y_3,.. ), \,(y_1,y_2,..)\,).$$

Then,
$$ \int \,
\sum_{a\in \{1,2\}} p_a\, k( \,( a,y_2,.. ), \,(y_1,y_2,..)\,)\,\, d \mu(y)= $$
$$ \int \mathcal{L}_{\log p} (g)(y_1,y_2,y_3,..)\, d\mu(y)=\int  k( \,(y_1,y_3,...)\,, \,( y_2,y_3,...)\,) d\mu(y).$$

Now, given $k$ consider the function
$$h(x_1,x_2,x_3,x_4,...)= k( \,(x_2,x_3,...)\,, \,( x_1,x_3,x_4,...)\,).$$

Then,
$$ \int \mathcal{L}_{\log p} (h)(x_1,x_2,x_3,..)\, d\mu(y)=\int  k( \,(x_2,x_3,...)\,, \,( x_1,x_3,...)\,) \, d \mu(x) .$$

For the Gibbs probability $\mu$ for $\log p$ is not true that for all $k$
$$\int  k( \,(x_2,x_3,...)\,, \,( x_1,x_3,...)\,) \, d \mu(x) = \int  k( \,(x_1,x_3,...)\,, \,( x_2,x_3,...)\,) d\mu(x).$$

\end{example}

\bigskip

\section{Noncommutative integration and quasi-invariant probabilities} \label{non}

\medskip

In non-commutative integration the transverse measures are designed to integrate transverse functions (see \cite{Con} or \cite{Kas}).

In the same way we can say that a function can be integrated by a measure resulting in a real number we can say that the role of a transverse measure  is to integrate transverse functions (producing a real number).

The main result here is Theorem \ref{mmeo} which describes a natural way to define a transverse measure from a modular function $\delta$ and a Haar system $(G\,,\hat{\nu})$.

As a motivation for the topic of this section consider a foliation of the two dimensional  torus where we denote each leave by $l$. This partition defines a grupoid with a quite complex structure. Each leave is a class on the associated equivalence relation. This motivation is explained with much more details in \cite{Con1}

We consider in each leave $l$ the intrinsic Lebesgue measure on the leave which will be denoted by $ \rho_l$.

A random operator $q$ is the association of a bounded operator $q(l)$ on $\mathcal{L}^2 (\rho_l)$ for each leave $l$. We will avoid to describe several technical assumptions which are necessary on the theory (see page 51 in \cite{Con1}).

The set of all random operators defines a von Neumann algebra under some natural definitions of the product, etc... (see Proposition 2 in page 52 in \cite{Con1})\,\,\,\,\,(*).

This setting is the formalism which is natural on {\bf noncommutative geometry} (see \cite{Con1}).

Important results on the topic are for instance   the characterization of when such von Neumann algebra is of type I, etc... (see page 53 in \cite{Con1}). There is a natural trace defined on this von Neumann algebra.

\bigskip

A more abstract formalism is the following:
consider a fixed grupoid $G$.
Given a transverse function $\lambda$ one can consider a natural operator $L_\lambda :\mathcal{F}^+(G)\to \mathcal{F}^+(G)$, which  satisfies
$$ f \to \lambda \,*\, f= L_\lambda (f)\,.$$

$L_\lambda $ acts on  $\mathcal{F}^+(G)$ and can be extended to a linear action on the von Neumann algebra $\mathcal{F}(G)$. This defines a  Hilbert module structure (see section 3.2 in \cite{Lan} or \cite{Kas}).

Given $\lambda$ we can also define the operator $R_\lambda :\mathcal{F}^+(G)\to \mathcal{F}^+(G)$ by
$$h(x,y)=  R_\lambda (f) (\gamma)\,=\, \int f(s, y) d\lambda^x ( d\,s ),$$
for any $(x,y).$

\begin{definition} \label{def9}
	
	Given two $G$-kernels $\lambda_1 $ and $\lambda_2$  we get a new $G$-kernel $\lambda_1 * \lambda_2$, called  the convolution of $\lambda_1 $ and $\lambda_2$, where given the function $f(x,y)$, we get the rule
	$$ (\lambda_1 * \lambda_2)\, (f)  \,=g \in  \mathcal{F} (G^0),$$
	given by
	$$g(y)=     \int\,(\,\,\int f(s, y)\,\,\lambda_2^x ( d\,s )\,) \,\,  \,\,\,\lambda_1^y   (dx)     .$$

	In the above $x \sim y \sim s$.

In other words $(\lambda_1 * \lambda_2)$ is such that for any $y$ we have
 \begin{equation} \label{kio}(\lambda_1 * \lambda_2)^y (dx)= \int \lambda_2^x ( d\,s\,) \,\,\,\lambda_1^y   (dx).
 \end{equation}
	
\end{definition}

Note that $$R_{\lambda_1 * \lambda_2} = R_{\lambda_1} \circ R_{\lambda_2}. $$

For a given fixed transverse function $\lambda$, for each class $[y]$ on the grupoid $G$, we get that
$R_\lambda $ defines an operator acting on functions $f(r,s)$, where $f: [y]\times [y] \to \mathbb{C},$ and
where $R_\lambda (f) =h $.

In this way, each transverse function $\lambda$ defines a random operator $q$, where $q([y])$ acts on $\mathcal{L}^2 (\lambda^y)$ via $R_\lambda$.

A transverse measure can be seen as an integrator of transverse functions or as an integrator of random operators (which are elements on the von Neumann algebra (*) we mention before).

\medskip

First we will present the basic definitions and results that we will need later on this section.

Remember that $\mathcal{E}^+$ is the set of transverse functions for the grupoid $G\subset X \times X$ associated to a certain equivalence relation $\sim$.

$\mathcal{F}^+( G) $ denotes the space of Borel measurable functions $f: G \to [0,\infty)$ (a real function of two variables $(a,b)$).

\begin{definition}\label{kern*func}

	Given a $G$ kernel $\nu$ and
	an integrable function $f\in \mathcal{F}_\nu (G)$ we can define two  functions on $G$:
	$$(x,y) \to(\nu  * f)(x,y) =\int f\,(x,s)\, \nu^y (ds)$$
	
	and
	
	$$ (x,y) \to (f * \nu)(x,y) =\int f\,(s,y)\, \nu^x (ds).$$

\end{definition}

Note that $\nu * 1=1$ if $\nu^y$ is a probability for all $y$.
Also note that $(f * \nu)(y,y)=  \nu  (f)\,\,(y)$ (see definition \ref{intkernel}).

\medskip
About (\ref{kio})
we observe that

\[  (\lambda_1 * \lambda_2)\, (f) = \lambda_1(f * \lambda_2) . \]
\medskip

A kind of analogy of the above concept of convolution (of kernels) with integral kernels is the following: given the kernels $K_1(s,x)$ and $K_2 (x,y)$ we define the kernel
$$ \hat{K} (s,y)  = \int K_2(s,x) \, K_1 (x,y)\, d x.$$ This is a kind of convolution of integral kernels.

This defines the operator
$$ f(x,y) \to g(y)= \int f(s,y) \hat{K}(s,y) d s= \int\,(\,  \int f(s,y)\,K_2(s,x) \,d\,s\,)\, \, K_1 (x,y)\, d x     .  $$

\begin{example} \label{delt} Given any kernel $\nu$ we have that $\mathfrak{d}\,* \,\nu=\nu$, where  $\mathfrak{d}$ is the delta kernel of Example \ref{del}.
	
	Indeed, for any $f \in \mathcal{F} (G)$
	$$  \int f \, (\,\mathfrak{d} * \nu\,)^y =\int \, \int f(s,y) \, \nu^x (d\,s) \, \mathfrak{d}^{\,y} (d\,x)=  \, \int f(s,y) \, \nu^y (d\,s) \,=\int f \, \nu^y$$
	
	In the same way for any $\nu$ we have that $\nu * \mathfrak{d}=\nu.$
	
\end{example}


\begin{example} \label{ess}  Given a fixed positive function $h(x,y)$ and a fixed kernel $\nu$, we get that the kernel $ \nu \, * (\,h\, \mathfrak{d}\,)$, where $\mathfrak{d}$  is the
	Dirac kernel, is such that given any $f(x,y)$,
	$$ (\nu \, * (\,h\,\mathfrak{d} \,))(f)(y)=     \int\,(\,\,\int f(s, y)\,h(s,x)\,\mathfrak{d}^x ( d\,s )\,) \,\,  \,\,\,\nu^y   (dx)     =$$
	$$  \int\, f(x, y)\,h(x,x)\, \,\,  \,\,\,\nu^y   (dx)     .$$
	
	Particularly, taking $h=1$, we get $\nu * \mathfrak{d} = \nu$.
	
\end{example}

\begin{example} \label{XY1} For the bigger than two equivalence relation of example \ref{XY} on
	$(S^1)^\mathbb{N}$, where $S^1$ is the unitary circle, the equivalence classes are of the form
	$ \{\,(a,x_2,x_3,...), \,a\in S^1\, \},$ where $ x_j\in S^1$, $j \geq 2$, is fixed.

	Given $x = (x_1,x_2,x_3,...)$ we define $\nu^x (da)$ the Lebesgue probability on $S^1$, which can be identified with $S^1 \times (x_2,x_3,..,x_n,...)$. This defines a transverse function where $G^0 = (S^1)^\mathbb{N}$. We call it the {\bf standard $\mathbf{XY}$ Haar system}.

	In this case given a function $f(x,y) = f(\,(x_1,x_2,x_3,..), (y_1,y_2,y_3,..)\,) $
	
	$$(\nu  * f)(x,y) =\int f(x,s)\, \nu^{x} (d\,s)=\int f\,(  \,(x_1,x_2,x_3,..), (s,x_2,x_3,..)\,)\, ds,$$
	where $s\in S^1$. Note that in the present example the information on $y$ was lost after convolution.
	
	Such $\nu$ is called in \cite{BCLMS} the a priori probability for the Ruelle operator. Results about Ruelle operators and Gibbs probabilities for such kind of $XY$ models appear in \cite{BCLMS}
	and \cite{LMMS2}.
	
\end{example}

\medskip

After Proposition \ref{kdf} we will present several properties of convolution of transverse function (we will need soon some of them).

\medskip

Note that if $\nu$ is transverse and $\lambda$ is a kernel, then $\nu * \lambda $ is transverse.

Remember that given a kernel $\lambda$ and a fixed $y$  the property  $\lambda^y(1) =1$
means $\int \lambda^y (dx)=1.$
\medskip

\begin{definition} \label{coco}
	A transverse measure $\Lambda$ over the modular function  $\delta(x,y)$, $\delta: G \to \mathbb{R},$ is a linear function  $\Lambda :\mathcal{E}^+\to \mathbb{R}^+$,  such that, for each kernel $\lambda$   which satisfies the property  $\lambda^y(1) =1$, for any $y$, if $\nu_1$ and $\nu_2$ are transverse functions such that $\nu_1 \,*\, (\delta \lambda)=\nu_2$, then,
	\begin{equation}
	\label{gogo}  \Lambda (\nu_1)= \Lambda (\nu_2).
	\end{equation}
	
\end{definition}

\medskip

A measure produces a real number from the integration of a classical function (which takes values on the real numbers),  and, on the other hand, the transverse measure produces a real number from a  transverse function $\nu$ (which takes values on measures).

\medskip

The assumptions on the above definition  are  necessary  (for technical reasons) when considering the abstract concept of integral of a transverse function by $\Lambda$ (as
is developed in \cite{Con}). We will show later that there is a more simple expression
providing the real values of such process of integration by $\Lambda$ which is related to quasi-invariant probabilities.

\medskip

If one consider the equivalence relation such that each point is related just to itself, any cocycle is constant equal
$1$ and the only kernel satisfying  $\lambda^y(1) =1$, for any $y$,  is the delta Dirac kernel $\mathfrak{d}.$
In this case if $\nu_1$ and $\nu_2$ are such that $\nu_1 \,*\, (\delta \lambda)=\nu_2$, then, $\nu_1=\nu_2$ (see Example \ref{delt}). Moreover, $\mathcal{E}^+$ is just the set of positive functions on $X$. Finally, we get that the associated transverse measure $\Lambda$  is just a linear function  $\Lambda :\mathcal{E}^+\to \mathbb{R}^+$

\medskip

\begin{example}
	Given a probability $\mu$ over $G^0$ we can define
	$$\Lambda(\nu)= \int \int \nu^y (dz) d \mu(y).$$
	
	Suppose that $\lambda$ satisfies $\lambda^x(1) =1$, for any $x$, and
	$$ \nu_1 \,*\, \lambda=\nu_2. $$
	
	Then, $\Lambda(\nu_1)= \Lambda(\nu_2).$ This means that $\Lambda$ is invariant by translation on the right side.
	
	\bigskip

	Indeed, note that,
	$$\Lambda(\nu_1)= \int \int \nu_1^y (dz) d \mu(y),$$
	and, moreover
	
	$$\Lambda(\nu_2)= \int \int \nu_2^y (dz) d \mu(y)=$$

	$$ \, \int [\,     \int\,(\,\,\int \,\,\lambda^x ( d\,s )\,) \,\, ) \,\,\,\nu_1^y   (dx) ) \,]\, d \mu(y)    =$$
	$$ \,      \int\,(\,\,\int  \,\,\,\nu_1^y   (dx) ) \,\, d \mu(y)    .$$

	Therefore, $\Lambda$ is a transverse measure of modulus $\delta=1$.
	
	In this way for each measure $\mu$ on $G^0$ we can associate a transverse measure of modulus $1$ by the rule $ \nu \to \Lambda(\nu)= \int \int \nu^y (dz) d \mu(y) \in \mathbb{R}.$
	
\end{example}

\bigskip

The condition
$$ \nu_1 \,*\, (\delta \lambda)=\nu_2$$

means for any $f$ we get
$$ \int f(x,y) \,(\,\nu_1 * \,(\delta \lambda \,))^y\,(d\,x)\,=  \int\,(\,\,\int f(s, y)\,[\,\delta (s,x) \lambda^x (ds)\,]\,)\, \,\,  \,\,\,\nu_1^y   (dx)     =$$
\begin{equation} \label{rr1} \int f(x,y) \,\,\nu_2^y (dx).
\end{equation}

\medskip

We define before (see (\ref{kwe})) the concept of quasi-invariant probability for a
given  modular function $\delta$, a grupoid $G$ and a fixed transverse function $\nu$.

 For reasons of notation we use a little bit variation of that definition. In this section we  say that  $M$ is quasi invariant probability  for $\delta$  and $\nu$ if for any $f(x,y)$

\begin{equation} \label{qui} \int \int f(y,x)\, \nu^y (x) d M( y) =\int \int f(x,y)  \delta(x,y)^{-1} \nu^y (x) d M(y). \end{equation}

\begin{proposition} \label{mmeo3} Given a modular function $\delta$, a grupoid $G$ and a fixed transverse function $\hat{\nu}$ denote by $M$ the quasi invariant probability for $\delta$.

Assume that $\int \hat{\nu}^y (dr) \neq 0$ for all $y$.

If $ \hat{\nu} * \lambda_1 = \hat{\nu} * \lambda_2 $, where $\lambda_1,\lambda_2$ are kernels, then,
$$ \int \delta^{-1} \lambda_1(1) \, \,d\, M= \int \delta^{-1} \lambda_2(1) \, \,d\, M. $$

This is equivalent to say that
$$\int \int \delta^{-1}(s,y)\lambda_1^y (ds) d M(y) =
\int \int \delta^{-1}(s,y)\lambda_2^y (ds) d M(y).$$

\end{proposition}

{\bf Proof:} By hypothesis $ g(y) = (\hat{\nu} * \lambda_1) (\delta^{-1})\, (y) = (\hat{\nu} * \lambda_2)(\delta^{-1}) (y). $

Then, we assume that (see (\ref{kio}))
$$\int g(y) \frac{1}{\int \hat{\nu}^y (dr) } d M(y) = \int \int \int \frac{1}{\int \hat{\nu}^y (dr) }  \delta^{-1} (s,y)\lambda_1^x (ds) \hat{\nu}^y (d x)\, d M(y)=$$
\begin{equation} \label{otite}\int \int \int \frac{1}{\int \hat{\nu}^y (dr) }  \delta^{-1}(s,y) \lambda_2^x (ds) \hat{\nu}^y (d x)\, d M(y).\end{equation}

Therefore,
$$ \int \int \delta^{-1}(s,y) \lambda_1 ^y (ds) d M(y)=$$
$$\int \int \int \frac{1}{\int \hat{\nu}^x (dr) } \delta(y,s) \lambda_1^y (ds)  \hat{\nu}^y (d x)\, d M(y) =$$
$$\int \int \int \frac{1}{\int \hat{\nu}^y (dr) } [\,\delta(y,x) \,\delta(x,s)\,\lambda_1^y (ds) \,] \hat{\nu}^y (d x)\, d M(y)=$$
$$\int \int \int \frac{1}{\int \hat{\nu}^y (dr) } [\,\delta(x,y)^{-1}  \,\delta(x,s)\,\lambda_1^y (ds) \,] \hat{\nu}^y (d x)\, d M(y)=$$
$$\int \int \int \frac{1}{\int \hat{\nu}^x (dr) } [\,\delta(y,x)^{-1}  \,\delta(y,s)\,\lambda_1^x (ds) \,]\,\delta^{-1} (x,y) \hat{\nu}^y (d x)\, d M(y)=$$
\begin{equation} \label{rsrs} \int \int \int  \frac{1}{\int \hat{\nu}^x (dr) } \delta(y,s) \lambda_1^x (ds)  \hat{\nu}^y (d x)\, d M(y)
\end{equation}

On the above from the fourth  to the fifth line we use the quasi-invariant expression (\ref{qui}) for $M$ taking
$$ f(y,x)= \int \frac{1}{\int \hat{\nu}^y (dr) } \,\delta(x,y)^{-1}  \,\delta(x,s)\,\lambda_1^y (ds) .$$

Note that if $\hat{\nu}$ is transverse $\int \hat{\nu}^x (dr)$ does not depends on $x$ on the class $[y]$.

Finally, from  the above equality (\ref{rsrs}) (and replacing $\lambda_1^x$ by $\lambda_2^x$)  it follows that
$$\int \int  \delta^{-1}(s,y)\lambda_1^y (ds)  \, d M(y) =$$
$$\int \int \int \frac{1}{\int \hat{\nu}^y (dr) } \delta (y,s) \lambda_1^y (ds)  \hat{\nu}^y (d x)\, d M(y) =$$
$$\int \int \int  \frac{1}{\int \hat{\nu}^y (dr) } \delta(y,s) \lambda_2^y (ds)  \hat{\nu}^y (d x)\, d M(y)=$$
$$\int \int  \delta^{-1}(s,y)\lambda_2^y (ds)  \, d M(y).$$

\qed

\bigskip

From now on we assume that $\int \hat{\nu}^y (dr) \neq 0$ for all $y$.

\medskip

\begin{theorem} \label{mmeo} Given a modular function $\delta$ and a Haar system $(G\,,\hat{\nu})$, suppose $M$ is quasi invariant for $\delta$.

	We define   $\Lambda$ on the following way: given a transverse function  $\hat{\nu} $ there exists a kernel $\rho$ such that $\nu= \hat{\nu} * \rho$ by proposition \ref{gos1}. We set
	\begin{equation}
	\label{gogoo5} \Lambda (\nu) = \int \int \delta(x,y)^{-1} \rho^y (dx) d M(y).
	\end{equation}

	Then, $\Lambda$ is well defined and it is a  transverse measure.

\end{theorem}

{\bf Proof:} $\Lambda$ is well defined by proposition \ref{mmeo3}.

We have to show that if $\lambda^x(1) =1$, for any $x$, and $\nu_1$ and $\nu_2$ are such that $\nu_1 \,*\, (\delta \lambda)=\nu_2$, then,
$\Lambda (\nu_1)= \Lambda (\nu_2).$

Suppose $\nu_1= \hat{\nu} * \lambda_1$, then, $\nu_2=   \hat{\nu} * \,(\,\lambda_1 * (\delta \,\lambda)\,)$.

Note that
$$\Lambda (\nu_1) = \int \int \delta(x,y)^{-1} \,\lambda_1^y (dx) d M(y).$$

On the other hand from (\ref{kio})
$$\Lambda (\nu_2) = \int \int \delta(s,y)^{-1} \,(\,\lambda_1 * (\delta \, \lambda\,)\, )^y (ds) d M(y)=$$
$$  \int \int \int \delta(s,y)^{-1} \,\delta(s,x)  \, \lambda^x (ds) \, \lambda_1^y (dx)\, d M(y)=$$
$$ \int \int \int \delta(x,y)^{-1}  \, \lambda^x (ds) \, \lambda_1^y (dx)\, d M(y)=$$
$$ \int \int \delta(x,y)^{-1} \,\, \,(\int\, \lambda^x (ds) \,)\,\,\, \lambda_1^y (dx)\, d M(y)=$$
$$ \int \int \delta(x,y)^{-1} \,\,\, \lambda_1^y (dx)\, d M(y)=\Lambda (\nu_1).$$
\qed

{\bf Remark:} The last proposition shows that given a quasi invariant probability $M$ - for a transverse function $\hat{\nu}$ and a cocycle $\delta$ - there is a natural way to define  a transverse measure $\Lambda$ (associated to a grupoid $G$ and a modular function $\delta$).

One can ask the question: given transverse measure $\Lambda$ (associated to a grupoid $G$ and a modular function $\delta$) is it possible to associate a probability on $G_0$? In the affirmative case, is this probability quasi invariant? We will elaborate on that.

\medskip

\begin{definition} Given a transverse measure $\Lambda$ for $\delta$ we can associate by Riesz Theorem  to a transverse function $\hat{\nu}$ a measure $M$ on $G^0$ by the rule: given a non-negative continuous function $h:G^0 \to \mathbb{R}$ we will consider the transverse function $h(x)\, \hat{\nu}^y (dx)$ and set
	
	$$ h \to \Lambda (h\, \hat{\nu} ) = \int h(x) d M(x).$$
	
	Such $M$ is a well defined measure (a bounded linear functional acting on continuous functions) and we denote such $M$ by $\Lambda_{\hat{\nu}}.$
	
	$\Lambda_{\hat{\nu}}$  means the rule $ h \to \Lambda (h\, \hat{\nu})= \Lambda_{\hat{\nu}} (h)$.

\end{definition}

\begin{proposition} \label{Conab}
	
	Given any transverse measure $\Lambda $ associated to the modular function $\delta $  and any transverse functions $\nu$ and $\nu'$ we have for any continuous $f$  that
	\[\Lambda_{\nu'}(\nu(\tilde \delta f))=\Lambda(\nu(\tilde \delta f) \nu' ) = \Lambda(\nu'(\tilde f) \nu)=\Lambda_{\nu}(\nu'(\tilde f)).  \]
	
\end{proposition}

{\bf Proof:} If $\lambda^y(1)=1 \,\, \forall y$, that is, $ \int 1 \lambda^y(ds) = 1\,\,\, \forall y$, then
$\Lambda(\nu*\delta \lambda) = \Lambda (\nu)$. If $g(x)=\lambda^x(1)=\int 1 \lambda^x(ds) \neq 1$, then we can write $\lambda'^x(ds)=\frac{1}{g(x)}\lambda^
x(ds)$, where $\lambda$ and $\lambda'$ are just kernels.  In this way $(\nu*\delta \lambda)= (g\nu)*\delta \lambda'$. Indeed, for $h(x,y)$,
\[\int h(x,y) \,(\nu*\delta \lambda)^y(dx) = \int h(s,y)\delta(s,x)\lambda^x(ds)\nu^y(dx) \]
\[= \int h(s,y)\delta(s,x)\lambda'^x(ds) g(x)\nu^y(dx)    =   \int h(x,y) ((g\nu)*\delta \lambda')^y(dx). \]

Denoting $\lambda(1) (x)=g(x)=\lambda^x(1)=\int 1 \lambda^x(ds)$, it follows that
\begin{equation}
\label{eq3}
\Lambda (\nu * \delta \lambda) = \Lambda (g\nu * \delta \lambda')=\Lambda (g\nu)=\Lambda_\nu(g)=\Lambda_\nu(\lambda(1)) =\Lambda_\nu(\int 1 \lambda^x(ds)) .
\end{equation}

From, (\ref{eq1}) if $\nu$ is a kernel and $f=f(x,y)$
\[ (\nu* f)(x,y) = \nu( \tilde f) (x),\]
and, from  Lemma \ref{gosy}, if $\lambda$ is a kernel and $\nu$ is a transverse function,  then, for any $f=f(x,y)$,
\[\lambda *(f\nu)=(\lambda*f)\nu.\]
It follows that, for transverse functions $\nu$ and $\nu'$, we get
\[ \nu * [(\delta \tilde f)\nu'] = [\nu*(\delta \tilde f)] \nu' = [\nu(\tilde \delta f)] \nu'.  \]

As a consequence
\[\Lambda_{\nu'}(\nu(\tilde \delta f)) = \Lambda([\nu(\tilde \delta f)]\nu') = \Lambda(\nu * [(\delta \tilde f)\nu']) =  \Lambda(\nu * \delta (\tilde f\nu'))= \]
\[  \Lambda_{\nu}((\tilde f \nu')(1)) =\Lambda_{\nu}(\int 1\cdot \tilde f (s,y) \nu'^y(ds)) = \Lambda_{\nu}(\nu'(\tilde f)) .  \]

Above we use  equation (\ref{eq3}) with $\lambda=\tilde f\nu'$.
\qed

\medskip

\begin{corollary} \label{simpl} If $\nu\in \mathcal{E}^+$, then for any $f$
	\begin{equation} \label{lopc} \Lambda ( \nu (\tilde{f})\,\nu\,)\, = \, \Lambda ( \nu (\tilde{\delta}\,f)\,\,\nu)
	\end{equation}
\end{corollary}

{\bf Proof:} Just take $\nu=\nu'$ on last Proposition.

\qed

\medskip

Among other things we are interested on  a modular function $\delta$,  a  transverse function $\hat{\nu}$ and a transverse measure $\Lambda$ (of modulo $\delta)$ such that   $M_{\Lambda,\hat{\nu}}=M$ is Gibbs for a Jacobian $J$. What conditions are required from $M$?

The main condition of the next theorem is related to the KMS condition of definition \ref{dede}

\medskip


\begin{proposition} \label{mme} Given  a transverse measure $\Lambda$ associated to the modular function $\delta$, and a {\bf transverse function} $\hat{\nu}$, consider the associated $M=\Lambda_{\hat{\nu}}$. Then, $M$ is quasi invariant for $\delta$. That is, $M$ satisfies for all $g$
	
	\begin{equation} \label{loip} \int \, \int g(s,x) \hat{\nu}^x (ds) d M(x)= \int \, \int g(x,s) \delta(x,s)\,\hat{\nu}^x (ds) d M (x).
	\end{equation}

\end{proposition}

{\bf Proof:}
 First we point out that (\ref{loip}) is consistent with (\ref{qui}) (we are just using different variables).

A transverse function $\hat{\nu}$ defines a function of $f\in \mathcal{F}(G) \to \mathcal{F}(G^0).$

The probability $M$ associated to $\hat{\nu}$ satisfies  for any continuous function $h(x)$, where $h:G^0 \to \mathbb{R}$ the rule
$$ h \to \Lambda (h\, \hat{\nu} ) = \int h(x) d M(x),$$
where $h(x)\, \hat{\nu}^y (dx)\in \mathcal{E}^+.  $

From proposition \ref{simpl}  we have that for the continuous function $f(s,x)=\tilde{g}(s,x)$, where $f:G \to \mathbb{R}$,
the expression
$$ \Lambda ( \hat{\nu}\,\,( g)\,\hat{\nu}\,)\,=\Lambda ( \hat{\nu}\,\,( \tilde{f})\,\hat{\nu}\,)\, = \, \Lambda ( \hat{\nu} (\delta^{-1}\,f)\,\,\hat{\nu}\,\,) =\Lambda ( \hat{\nu} (\delta^{-1}\,\tilde{g})\,\,\hat{\nu}\,\,)$$

For a given  function $g(s,x)$ it follows from the above  that
$$ \Lambda ( \hat{\nu}\,\,( g)\,\hat{\nu}\,)\,= \int \hat{\nu}\,\,( g)(x)\, d M(x)=\int\, \int g(s,x) \hat{\nu}^x (ds) \, d M(x).$$

On the other hand
$$\Lambda ( \hat{\nu} (\delta^{-1}\,\tilde{g})\,\,\hat{\nu}\,\,)=\int \hat{\nu} (\delta^{-1}\,\tilde{g}) (x) d\,M(x)=\int g(x,s) \delta^{-1} (s,x)\,\hat{\nu}^x (ds) \, d M(x).$$

\qed

\medskip

\begin{proposition} \label{mmeo1} Given a modular function $\delta$, a grupoid $G$, a transverse measure $\Lambda$ and a transverse function $\hat{\nu}$, suppose for any $\nu$, such that $\nu= \hat{\nu} * \rho$, we
	have that
	$$ \Lambda(\nu) = \Lambda (\hat{\nu} * \rho) \,= \int \int \delta(s,x)^{-1} \rho^x (ds) d \mu(x)= \int \delta^{-1} \rho(1) \, d\, \mu.  $$

	Then,
	$\mu= \Lambda_{\hat{\nu}}$.
	
\end{proposition}

{\bf Proof:} Given $f\in \mathcal{F}( G_0)$ consider $\lambda$ the kernel such that $ \lambda^x (d s)= f(x) \delta_x (ds)$, where $\delta_x$ is the Delta Dirac on $x$.

Then, using the fact that $\delta(x,x)=0$ we get that the kernel  $f(x) \hat{\nu}^y( dx)$ is equal to $\hat{\nu} * \delta \, \lambda$.

Then, taking $\rho= \delta\, \lambda$ on the above expression we get
$$ \Lambda (\, f \hat{\nu}  )=  \Lambda (\hat{\nu} * \,(\delta \,\lambda)\,) =$$
$$\int \delta^{-1} \rho(1)   d \mu= \int \delta^{-1}\,( \delta \lambda) (1)  d \mu =\int \lambda(1) d \mu =\int f(x) d \mu(x) . $$

Therefore, $\Lambda_{\hat{\nu}}= \mu.$

\qed

\medskip

Now we present a general procedure to get transverse measures.

\begin{proposition} \label{kdf}  For a fixed modular function $\delta$ we can associate to any  given  probability $\mu$ over $G^0$ a transverse measure  $\Lambda$ by the rule
	\begin{equation}\nu \to  \label{tri} \Lambda (\nu)=  \,\int \int \delta(s,x)^{-1} \,  \nu^x (ds) \,\, \, d \mu (x).
	\end{equation}

\end{proposition}

{\bf Proof:}

Consider $\nu\,'\in \mathcal{E}^+$ and $\lambda $, such that, $\int \lambda^r( ds )=1$, for all $r$, and moreover that $\nu\,' = \nu * (\delta \lambda)$.

We will write
$$ (\nu * \, \delta \lambda)\, (\delta^{-1} )= \int \,\int \delta^{-1} (s,x) \, \delta(s,r)\, \lambda^r (ds) \nu^x (dr)$$
which is a function of $x$

Then
$$ \Lambda (\nu\,')=  \,\int \int \delta(s,x)^{-1} \,  \nu\,'^x (ds) \,\, \, d \mu (x)=\int \nu\,'(\delta^{-1})(x)d \mu(x)=$$
$$=\int  \Big(\nu * (\delta \lambda)\Big)(\delta^{-1})(x)d \mu(x) =\int \int \int \delta(s,x)^{-1}\delta(s,r)\lambda^r(ds)\nu^x(dr) d \mu(x)= $$
$$=\int \int \int \delta(r,x)^{-1}\lambda^r(ds)\nu^x(dr) d \mu(x)=\int \int  \delta(r,x)^{-1}\nu^x(dr) d \mu(x)=\Lambda(\nu). $$

\qed

This last transverse measure is defined in a quite different  way that the one described on Theorem \ref{mmeo}.

\medskip

Now we will present some general properties of convolution of transverse functions.

\medskip

\begin{lemma} \label{gos} Suppose $\nu \in \mathcal{E}^+ $ is a transverse function, $\nu_0$ a kernel, and $g \in \mathcal{F}^+ (G)$ is such that $\int g(s,x) \nu_0^y (dx) =1$, for all $s,y$. Then,   $ \nu_0\, * \,(g \, \nu)\,=  \nu\,$, where
	$g \, \nu$ is a kernel.
\end{lemma}

\textbf{Remark} The condition $\int g(s,x) \nu_0^y (dx) =1$, for all $s,y$ means $(\nu_0 *g)(s,y) =1$ for all $s,y$, that is   $\nu_0 *g \equiv 1$ (See lemma 3 below).

\bigskip

{\bf Proof:}

$$ z(y)=\int \,f(s,y)\,  \,\nu^y (ds) = $$
$$\int f(s, y)\,\,[\,\int g (s,x)  \nu_0^y   (dx)]\, \nu^x (ds)\,\,\, \,\,  \,\,\,=$$
$$\int\,\,\,\int f(s,y)\,\,[\,g (s,x) \nu^x (ds)\,]\,\, \,\,  \,\,\,\nu_0^y   (dx)=$$
$$    \int f(s,y)\,   (g\,\nu)^x  (ds) \nu_0^y (dx)   =  \int f(s,y) \,(\,\nu_0 * \,(g \, \nu))^y (ds)  .$$

\qed

\medskip

We say that the kernel $\nu$ is fidel if $\int \nu_0^y (ds)\neq 0$ for all $y$.

\begin{proposition} \label{gos1}
	For a fixed transverse function $\nu_0$  we have that for each given  transverse function $\nu$ there exists a kernel $\lambda,$ such that,
	$\nu_0 * \lambda = \nu$.
\end{proposition}

{\bf Proof:} Given the kernel $\nu_0$ take $g_0(s)= \frac{1}{ \int 1 \nu_0^s (dr)} \geq 0$. Note that $g_0(v)$ is constant for $v\in [s]$. Then   $\nu_0(g)=1$, that is, for each $s$ we get that $\int g_0(s) \nu_0^s(dx)=1$.

We  can take $\lambda=\,g_0 \,\nu  $ as a solution. Indeed, in a similar way as last lemma we get

$$ z(y)=\int \,f(s,y)\,  \,\nu^y (ds) = $$
$$\int f(s, y)\,\,[\,\int g_0 (s)  \nu_0^s   (dx)]\, \nu^x (ds)\,\,\, \,\,  \,\,\,=$$
$$\int\,\,\,\int f(s,y)\,\,[\,g_0 (s) \nu^x (ds)\,]\,\, \,\,  \,\,\,\nu_0^y   (dx)=$$
$$    \int f(s,y)\,   (g_0\,\nu)^x  (ds) \nu_0^y (dx)   =  \int f(s,y) \,(\,\nu_0 * \,(g_0 \, \nu))^y (ds)  =$$
$$
\int f(s,y) \,(\,\nu_0 * \,\lambda)^y (ds)  .$$

\qed

The next Lemma is just a more general form of Lemma \ref{gos}.

\begin{lemma} \label{gosy} Suppose $\nu \in \mathcal{E}^+$, $g \in \mathcal{F}^+ (G)$ and $\lambda $ a kernel, then   $ \lambda\, * \,(g \, \nu)\,= ( \lambda \,* \,g )\,\,\, \nu\,$, where
	$g \, \nu$ is a kernel and $\lambda * g $ is  a function.
	
\end{lemma}

{\bf Proof:}

Given $f\in \mathcal{F}(G)$  we get
$$ (\,\lambda * \,(g \, \nu))(f)(y)= \int f(x,y) \,(\,\lambda * \,(g \, \nu))^y(dx) $$
$$=  \int\int f(s, y)[(g \nu)^x (ds)]\,\lambda^y (dx)=\int\int f(s, y)[g (s,x) \nu^x (ds)]\,\lambda^y   (dx)     .$$

On the other hand

$$[( \lambda \,* \,g )\,\,\, \nu](f)(y) = \int f(s,y) [(\lambda * g)\, \nu]^y (ds)= \int f(s,y) [(\lambda * g)(s,y)] \nu^y (ds)$$
$$= \int f(s,y) [\int g (s,x) \lambda^y (d x)] \nu^y (ds) = \int \int f(s,y)  g (s,x) \nu^x (ds) \lambda^y (d x) .$$

\qed

\medskip

\begin{proposition} \label{fifi} Suppose $\nu$ and $\lambda$ are transverse. Given $f\in \mathcal{F}^+ (G)$, we have that
	$$ \lambda\,(\nu * f)   \,=\,  \nu\,(\lambda * \tilde{f}).$$
	
\end{proposition}

{\bf Proof:}

Indeed, by definition \ref{kern*func}, $(\nu * f)(x,y)=g(x,y)=\int f(x,s)\nu^y(ds)$, and by definition \ref{intkernel} $$\lambda\,(\nu * f)(y)=\lambda(g)(y)=\int g(x,y)\lambda^y(dx)=\int \int f(x,s)\nu^y(ds) \lambda^y(dx).$$
By the same arguments $(\lambda * \tilde{f})(x,y)=h(x,y)=\int \tilde f(x,s)\lambda^y(ds)$, and
$$\nu\,(\lambda * \tilde f)(y)=\nu(h)(y)=\int h(x,y)\nu^y(dx)=\int \int \tilde f(x,s)\lambda^y(ds) \nu^y(dx)=$$
$$=\int \int  f(s,x)\lambda^y(ds) \nu^y(dx)= \lambda\,(\nu * f)(y), $$
if we exchange the coordinates $x$ and $s$.

Note that in the case $f\in \mathcal{F} (G^0)$ we denote $f(x,s)=f(x)$.
In the same way $$ \lambda\, (\nu * f)  \,=\,\nu \,(\lambda * \tilde{f}) $$
in the following sense:
$$ \int \int f(x)\nu^y(ds) \lambda^y(dx)=\int \int  f(s)\lambda^y(ds) \nu^y(dx)  .$$

\qed

\medskip

\medskip

\medskip

\section{ $C^*$-Algebras derived from Haar Systems} \label{C2}

\medskip

In this section the functions $f:G \to \mathbb{R}$ will be required to be continuous (not just measurable).

An important issue here is to have suitable hypotheses in such way that the indicator of the diagonal $\mathfrak{1}$ belongs to  the
underlying space we consider. On von Neumann algebras the unit  is just measurable and not continuous (this is good enough).
We want to consider another setting (certain $ C^*$-algebras associated to Haar Systems) where the unit will be required to be a continuous function. In general terms, given a groupoid $G\subset \Omega\times \Omega$, as we will see, we will need another topology on the set $G$ for the $C^*$-Algebra formalism and for defining KMS states.

We will begin with some more examples. The issue here is to set a certain appropriate topology.

\begin{example} \label{wun}
	
	For $n \in \mathbb{N}$ we define the partition $\eta_n  $ over $\overrightarrow{\Omega}=\{1,2,...,d\}^\mathbb{N}$, $d \geq 2$, such that two elements $x\in \overrightarrow{\Omega}$ and $y\in\overrightarrow{\Omega}$ are on the same element of the partition, if and only if, $ x_j=y_j, $ for all $j> n .$  This defines an equivalence relation denoted by $R_n$.
	
\end{example}

\begin{example} \label{wui0}

	We define a partition $\eta  $ over $\overrightarrow{\Omega}$, such that two elements $x\in \overrightarrow{\Omega}$ and $y\in\overrightarrow{\Omega}$ are on the same element of the partition, if and only if, there exists an $n$ such that $ x_j=y_j, $ for all $j> n .$  This defines an equivalence relation denoted by $R_\infty$.
	
\end{example}

\begin{example} \label{winf} For each fixed $n\in \mathbb{Z}$ consider the equivalence relation on $\hat{\Omega}$: $x \sim y$ if
	
	$$y =(..., y_{-n}, ... ,y_{-2},y_{-1}\,|\, y_{0},y_{1},....,y_{n},...)\,\,$$ is such that $ x_j=y_j $ for all $j\leq n ,$ where  $\hat{\Omega}=\overleftarrow{\Omega}\times  \overrightarrow{\Omega}.$
	
	This defines a  groupoid.

\end{example}

\begin{example}
	
	Recall that by definition the unstable set of the point $x\in\hat{\Omega}$ is the set
	$$W^u (x)=\{y \in \hat{\Omega}\,,\,\,\text{such that}\,\,\, \lim_{n\to \infty} d( \hat{\sigma}^{-n}(x), d( \hat{\sigma}^{-n}(y)\,)\,=0\,\}$$
	
	One can show that the unstable manifold of $x \in \hat{\Omega}$ is the set
	$$W^u (x)=\{y =(..., y_{-n}, ... ,y_{-2},y_{-1}\,|\, y_{0},y_{1},....,y_{n},...)\,|\,\text{there exists}\,$$
	$$ k\in \mathbb{Z},\,\, \text{such that}\,\, x_j=y_j, \text{for all}\,\,j\leq k  \}.$$

	If we  denote by $G_u$ the groupoid   defined by the above relation, then, $x \sim y$, if and only if $y \in W^u (x).$
	

\end{example}

\begin{definition} \label{pro}
	Given the equivalence relation  $R$, when the quotient $\hat{\Omega}/R$ (or, $\overrightarrow{\Omega}/R$) is  Hausdorff and locally compact  we say that $R$ is a proper equivalence.
	
\end{definition}

\medskip

For more details about proper equivalence see section 2.6 in \cite{Ren2}.

\medskip

On the set $ X=\overrightarrow{\Omega}$, if we denote $x=(x_1,x_2,..x_n,..) $,  the family $U_x(m)=\{y\in \overrightarrow{\Omega}, $ such that, $ y_1=x_1,y_2=x_2,...,y_m=x_m\}$, $m=1,2,...$, is a fundamental set of open neighbourhoods on $\Omega$.

Considering the relations $R_m$ and $R_\infty$ we get the corresponding groupoids
$$ G_1 \subset G_2\subset ... \subset  G_m\subset ... \subset G_\infty \subset \overrightarrow{\Omega} \times \overrightarrow{\Omega} =X \times X.$$

\bigskip

The equivalence relation described in example \ref{wui0} (and also \ref{win}) is not proper if we consider the product topology on $\overrightarrow{\Omega}$ (respectively on $\hat{\Omega}$). The equivalence relation described in example \ref{wun} (and also \ref{winf}) is proper if we consider the product topology on $\overrightarrow{\Omega}$ (respectively on $\hat{\Omega}$) (see \cite{EL1}).

\medskip

We consider over $G_n$ the quotient topology.

\begin{lemma} \label{lle1} Given $X= \overrightarrow{\Omega}$,
	for each $n$ the map defined by the canonical projection $ X \to G_n$ is open.
	
\end{lemma}

{\bf Proof}:  Given an open set $U\subset X$ take $V=\{y\in X\,|\, $ there exists $x\in X$, satisfying $ y\sim x$ for the relation $R_n\}$. We will show that $V$ is open.

Consider $y\in V$, $y\in U$, such that, $ y\sim x$ for the relation $R_n$. There exists $m>n$, such that, $U_x(m) \subset U$. Then,  $U_y(m) \subset V$. Indeed, if $z \in U_y(m)$, take $z' \in X$, such that $z_j ' =x_j,$ when $1\leq j\leq m$, and $z_j ' =z_j$ , when $j>m$.

Then, $z' \sim z$ for the relation $R_n$. But, as $ z \in U_y(m)$, this implies that $z_j=y_j$, when $1\leq j\leq m$, and
$y \sim x$, for $R_n$, implies that $y_j=x_j$, when $j>n$. Then, $z_j ' =x_j$, if $1 \leq j\leq m$. Therefore, $z' \in U_x(m) \subset U.$

\qed

\begin{lemma} \label{lle2} Given $X= \overrightarrow{\Omega}$,
	for each $n$ the map defined by the canonical projection $ X \to G_\infty$ is open.
	
\end{lemma}

{\bf Proof}:  Given an open set $U\subset X$ take $V=\{y\in X\,|\, $ there exists $x\in X$, satisfying $ y\sim x$ for the relation $R_n\}$  and   $V_\infty =\{y\in X\,|\, $ there exists $x\in X$, satisfying $ y\sim x$ for the relation $R_\infty\}$. Then, $V_\infty = \cup_{n=1}^\infty\, V_n$ is open.

\qed

\begin{lemma} \label{lle3}  Given $X= \overrightarrow{\Omega}$,
	for each $n=1,2...,n,...$, the set $G_n$ is Hausdorff.
	
\end{lemma}

{\bf Proof}: Given a fixed $n$, and  $x,y\in X$, such that $x$ and $y$ are not related by $R_n$, then, there exists $m>n$ such that $x_m \neq y_m$. From this follows that no element of $U_x(m)$ is equivalent by $R_n$ to an element of $U_y(m).$ By lemma \ref{lle1} it follows that $G_n$ is Hausdorff.

\qed

\begin{lemma} \label{lle4}  Given $X= \overrightarrow{\Omega}$ the set $G_\infty$ is not Hausdorff.
	
\end{lemma}

{\bf Proof}:  If $x_m = (\underbrace{1,1,...,1}_m,d,d,d...)$, then $\lim_{n \to \infty} x_n = (1,1,1...,1,...)$ and
$(\underbrace{1,1,...,1}_m,d,d,d...)\sim  (d,d,d,...,d,...)$, for the relation $R_\infty$. Note, however, that
$(1,1,1...,1,...)$ is not in the class $(d,d,d,...,d,...)$ for the relation $R_\infty$.

\qed

\begin{lemma} \label{lle5}  Given $X= \overrightarrow{\Omega}$ denote by $D$ the diagonal set on $X \times X$. Then, $D$ is open on $G_n$  for any $n$, where we consider on $D$ the topology induced by $X \times X$.

\end{lemma}

{\bf Proof}:  Given $x\in X$, we have that $U_x(n) \times U_x (n)$ is an open set of $X\times X$ which contains $(x,x)$.

Consider $y,z\in U_x (n)$ such that $y$ and $z$ are related by $R_n$. Then,
$y_j  =x_j= z_j,$ when $1\leq j\leq n$, and $y_j  =z_j$ , when $j>n$. Therefore, $y=z$.

From this we get that
$$ U_x(n) \times U_x (n)\cap G_n \subset D$$

\qed

\medskip

\begin{definition}
	An equivalence relation $R$ on a compact Hausdorff space X is said to be
	approximately proper if there exists an increasing sequence of proper equivalence relations
	$R_n$, $n \in \mathbb{R}$, such that $R=\cup_n R_n$, $n \in \mathbb{N}.$ This in the sense that if $x\sim_R y$, then there exists an $n$ such that $x\sim_{R_n} y$.
\end{definition}

\begin{example}
	Consider the equivalence relation $R_\infty$ of example \ref{wui0} and $R_n$ the one of example \ref{wun}. For each $n$ the equivalence relation $R^n$ is proper.
	
	Then, $R_\infty=\cup_n R_n$, $n \in \mathbb{N}$ is approximately proper  (see \cite{EL1}).

\end{example}

\begin{definition} \label{rtu}
	
	Consider a fixed set $K$, a sequence of subsets $W_0 \subset W_1 \subset W_2 \subset...\subset\, W_n \subset ...\subset K$ and a topology $\mathcal{W}_n$ for each set $W_n\subset K.$
	
	By the direct inductive limit
	$$t-\lim_{n \to \infty}\mathcal{W}_n= \mathcal{K}$$
	we understand the set $K$ endowed with the
	largest topology $\mathcal{K}$ turning the identity inclusions $W_n \to K$ into continuous maps.
	
	The topology of $t-\lim_{n \to \infty}\mathcal{W}_n= \mathcal{K}$ can be easily described: it consists of all subsets $U \subset K$
	whose intersection $U \cap W_n $  is in $\mathcal{W}_n$  for all $n$.

\end{definition}

For more details about  the inductive limit (see section 2.6 in \cite{Ren2}).

\medskip
In the case $W_n=G_n$ we consider as $\mathcal{W}_n$ the product topology.

\medskip

\begin{lemma} \label{lle6}  Given $X= \overrightarrow{\Omega}$ if we consider over $K=G_\infty$ the inductive limit topology defined by the sequence of the $G_n \subset X \times X$, then, the indicator function $\mathfrak{1}$  on the diagonal is continuous.

\end{lemma}

{\bf Proof}:  By lemma \ref{lle5} the diagonal $D$ is an open set.

Moreover, $(G_\infty - D)\,\cap G_n = (\, (X \times X) - D) \,\cap \,G_\infty)\, \cap G_n= (  \,(X \times X \,)\,-\,D ) \cap G_n$
is open on $G_n$ for all $n$. Then, $(G_\infty - D)$ is open on $G_\infty.$

\qed

\medskip

{\bf Remark:} Note that on $G_\infty$ we have that $D$ is not open on the induced topology by $X \times X$. Indeed,
consider $a=(1,1,1,..,1,..)$ and $b_m=  (\underbrace{1,1,...,1}_{m-1},d,1,1,1,..1,...)$. Then,
$\lim_{m \to \infty} (a,b_m)=(a,a)\in D$, and $(a,b_m)\in G_\infty$ but $(a,b_m)$ is not on $D$, for all $m$.

\bigskip

\begin{example} \label{uni}

	In the above definition \ref{rtu} consider $ W_n=G_n\subset \overleftarrow{\Omega}\times  \overrightarrow{\Omega} $, $n \in \mathbb{N}$, which is the groupoid associated to the equivalence relation $R_n$ (see example \ref{wun}).  Then, $\cup_n W_n=K=G\subset \overleftarrow{\Omega}\times  \overrightarrow{\Omega} $, where $G$ is the groupoid associated to the equivalence relation $R^\infty$. Consider on $W_n$ the topology $\mathcal{W}_n$ induced by the product topology on
	$ \overrightarrow{\Omega} \times  \overrightarrow{\Omega}$.
	
	For a fixed $x$ the set $U= \{y\,|\, x_j=y_j\,$ for all $j\leq n\}\cap G_n$ is open on $G_n$, that is, an element on $\mathcal{W}_n$.

	Note that $G_n \cap ( U \times U)$ is a subset of the diagonal.
	
	Points of the form
	$$(\, (x_1,x_2,...,x_n, z_{n+1},z_{n+2},...),  (x_1,x_2,...,x_n, z_{n+1},z_{n+2},...)\,)$$
	are on this intersection.
	
	Then, the diagonal $\{(y,y)\,,\,y \in \overrightarrow{\Omega}\}$ is an open set in the inductive limit topology $\mathcal{K}$ over $G$
	
	From this follows that the indicator function of the diagonal, that is, $I_\Delta$, where $\Delta = \{(x,x) |\, x \in \overrightarrow{\Omega}\}$, is a continuous function.

\end{example}

\begin{example} \label{wui}

	Consider the partition $\eta_n  $, $n \in \mathbb{Z}$, over $\hat{\Omega}$ of Example \ref{winf}, $W_n=G_n$, for all $n$, and  $K=G_u $.

	We consider the topology $\mathcal{W}_n$ over $G_n$ induced by the product topology. In this way $A \in \mathcal{W}_n$ if
	$$ A= B\cap G_n,$$
	
	where  $B$ is an open set on the product topology for $\hat{\Omega}\times \hat{\Omega} .$
	
	In this way $A$ is open on $t-\lim_{n \to -\,\infty}\mathcal{W}_{n}= \mathcal{K}$ if for all $n$ we have that
	$$ A \cap X_n \in \mathcal{X}_n.$$

	Denote by $D$ the diagonal on $\hat{\Omega}\times \hat{\Omega}$
	and consider the indicator function $I_D: \hat{\Omega}\times \hat{\Omega}\to \mathbb{R}$.
	
	The function $I_D$ is continuous over   the inductive limit topology $\mathcal{K}$ over $K=G_u$.
	
\end{example}

\medskip

Here $G^0$ will be the set $\hat{\Omega} =\overleftarrow{\Omega}\times  \overrightarrow{\Omega}$. We will denote by $G$ a general groupoid obtained by an equivalence relation $R$.

The measures we consider on this section are defined over  the sigma-algebra generated by the inductive limit topology.
\medskip

\begin{definition}
	Given a Haar system $(G,\nu)$, where $G^0=\hat{\Omega}= \overleftarrow{\Omega}\times  \overrightarrow{\Omega}$ is equipped with the  inductive limit topology, considering  two continuous  functions with compact support $f,g \in  C_C (G)$, we  define $(f \,  \underset{\nu}{*}\, g)=h$ in such way that for any $(x,y)\in G$
	$$(f \,  \underset{\nu}{*}\, g)(x,y)=
	\int \,g(x,s ) \,f (s,y) \,\,\nu^{y} (ds) =h(x,y).$$
	
	The closure of the operators of left multiplication by elements of $C_C(G)$, $\{L_f:f\in C_C(G)\}\subseteq B(L^2(G,\nu))$, with respect to the norm topology is called the reduced C*-algebra associated to $(G,\nu)$ and denoted by $C_r^*(G,\nu)$.
\end{definition}

{\bf Remark:} There is another definition of a C*-algebra associated to $(G,\nu)$ called the full C*-algebra. For a certain class of groupoid, namely the amenable groupoids, the full and reduced C*-algebras coincide. See \cite{ADR} for more details.

\bigskip

As usual function of the form $f(x,x)$ are identified with functions $f: G^0 \to \mathbb{C}$ of the form $f(x).$

The collection of these functions  is commutative sub-algebra of
the $C^*$-algebra $C^*_r(G,\nu)$.

We denote by $\mathfrak{1}$ the indicator function of the diagonal on $G^0 \times G^0$.
Then, $\mathfrak{1}$ is the neutral element for the product $\underset{\nu}{*}$ operation. Note that
$\mathfrak{1}$ is continuous according to example \ref{wui}.

In the case there exist a neutral multiplicative element we say the $C^*$-Algebra is unital.

Similar properties to the von Neumann setting can also be obtained.

We can define in analogous way to definition \ref{KMSdef} the concept of $C^*$-dynamical state (which requires an unit $\mathfrak{1}$) and the concept of KMS state for a continuous modular function $\delta$.

\medskip

General references on  the $C^*$-algebra setting are \cite{Ren0},  \cite{Ren2}, \cite{Exel1}, \cite{Exel2}, \cite{Exel3}, \cite{EL1}, \cite{EL2}, \cite{ManS}, \cite{PutJ}, \cite{KumRen} and \cite{AHR} .

\medskip

\medskip
\section{Examples of quasi-stationary probabilities} \label{exaquasi}

On this section we will present several examples of measured groupoids, modular functions and  the associated quasi-stationary probability (KMS probability).

\begin{example} \label{out}
 Considering the example \ref{inv1} we get that each $a   \in \{1,2,...,d\}^\mathbb{N}= \overleftarrow{\Omega}$
 defines a class of equivalence
 $$a\times \,|\, \overrightarrow{\Omega}=a  \times \{1,2,...,d\}^\mathbb{N}=  (..., a_{-n}, ... ,a_{-2},a_{-1}) \times \,|\, \{1,2,...,d\}^\mathbb{N}.$$

 On next theorem we will denote  by $G$ such  groupoid.

Given a Haar system $\nu$ over  such $G\subset \hat{\Omega} \times \hat{\Omega}$,
note that if $z_1=<a|b_1>$ and $z_2=<a|b_2>$, then $\nu^{z_1}= \nu^{z_2}.$ In this way it is natural do index the Haar system by $\nu^a$, where $a \in \overleftarrow{\Omega}.$ In other words, we have
\begin{equation} \label{puo} \nu^{ <a|b> }  (d\, <{a}|\tilde{b}>)=  \nu^{ a }  (d\, \tilde{b}).
\end{equation}

 Consider $V: G \to \mathbb{R}$,
 $m $ a probability over    $\overleftarrow{\Omega}$ and the modular function $\delta(x,y)= \frac{e^{V(x)}}{e^{V(y)}},$ where $(x,y)\in G$.

 Finally we denote by $\mu_{m,\nu,V}$ the  probability on $G^0=\hat{\Omega}$, such that, for any function $g:\hat{\Omega} \to \mathbb{R}$ and $y=<a\,|\,b>$
$$ \int g(y) d \mu_{m,\nu,V} (y) =\,\,  \int_{\overleftarrow{\Omega}}\,\,\,(\int_{\overrightarrow{\Omega}}\,\,g(<a|b>)\, e^{V(<a|b>)}\,d\nu^a(d\,b)  \,)\,\,d\, m(d a) .$$

Note that $\hat{\nu}= e^V \nu$ is a $G$-kernel but maybe not transverse.
The next theorem will provide a large class of examples of quasi-invariant probabilities for such groupoid $G$.
\end{example}

\bigskip
 \begin{theorem} \label{gou} Consider a Haar System $(G,\nu)$ for the groupoid of example \ref{out}. Then, given $m$, $V$, using the notation above we get that $M=  \mu_{m,\nu,V}$ is quasi-invariant for the modular function $\delta(x,y)=\frac{e^{V(x)}}{e^{V(y)}}$.

\end{theorem}

{\bf Proof:}


From (\ref{kwe1}) we just have to prove that for any $f: G \to \mathbb{R}$

\begin{equation} \label{kwe8} \int  \int f(x,y) e^{ V(x)} \nu^y (dx) d   \mu_{m,\nu,V}(dy)= \int  \int f(y,x) e^{V(x)}\nu^y (dx) d \mu_{m,\nu,V}(dy) .
\end{equation}

We denote $y=<a|b>$ and   $x=<\tilde{a}|\tilde{b}>.$ Note that if $y\sim x$, then $a=\tilde{a}.$

Note that, from (\ref{puo})
$$   \int\,\,(\int f\,(x,y)\,e^{V(x)} \nu^y (dx)\,\,)\,\, d \mu_{m,\nu,V}(dy) =$$
{\tiny $$   \int\,\,(\int f\,( <\tilde{a}|\tilde{b}>,\,<a|b>)\,e^{ V( <\tilde{a}|\tilde{b}> )}\nu^{\,<a|b>} (d  <\tilde{a}|\tilde{b}>\,)\,\, d \mu_{m,\nu,V}(d\,<a|b>) =$$
$$   \int_{\overleftarrow{\Omega}}\,[\,\,\int_{\overrightarrow{\Omega}}\,(\int\, \int \, f\,(<\tilde{a}|\tilde{b}> ,\,<a|b>)\,  \, e^{V( <\tilde{a}|\tilde{b}> )}\nu^{\,<a|b>} (d<\tilde{a}|\tilde{b}>)\,\, \,  e^{V(<a|b>)}\,d\nu^a(d\,b)  \,)\,\,]\,d\, m(d a)  =$$
$$   \int_{\overleftarrow{\Omega}}\,[\,\,\int_{\overrightarrow{\Omega}}\,(\int\, \int \, f\,(<a|\tilde{b}> ,\,<a|b>)\,  \,e^{ V( <a|\tilde{b}> )}\, \,  e^{V(<a|b>)}\,\nu^{\,a} (d\, \tilde{b})\,d\nu^a(d\,b)  \,)\,\,]\,d\, m(d a)  .$$}

In the above expression we can exchange the variables $b$ and  $\tilde{b}$, and,  finally, as $a = \tilde{a}$,  we get
{\tiny $$   \int_{\overleftarrow{\Omega}}\,[\,\,\int_{\overrightarrow{\Omega}}\,(\int\, \int \, f\,(<a|b> ,\,<a|\tilde{b}>)\,  \,e^{ V( <a|b> )}\, \,  e^{V(<a|\tilde{b}>)}\,\nu^{\,a} (d\, b)\,d\nu^a(d\,\tilde{b})  \,)\,\,]\,d\, m(d a)  =$$
$$   \int_{\overleftarrow{\Omega}}\,[\,\,\int_{\overrightarrow{\Omega}}\,(\int\, \int \, f\,(<a|b> ,\,<a|\tilde{b}>)\,  e^{V(<a|\tilde{b}>)}\, d\nu^a(d\,\tilde{b})\, \,\, \,e^{ V( <a|b> )}\,\nu^{\,a} (d\, b)\,  \,)\,\,]\,d\, m(d a)  =$$
$$   \,(\int\, \int \, f\,(y ,\,x)\,   e^{V(x)}\,\, d\nu^{y}(dx)\, \,d \mu_{m,\nu,V}(dy).$$}

This shows the claim.

\qed

\bigskip

\begin{example} \label{go} Consider $G$ associated to the equivalence relation given by the unstable manifolds for $\hat{\sigma}$ acting on $\hat{\Omega}$ (see example \ref{win}).
Let's fix for good a certain $x_0 \in \overrightarrow{\Omega}$.
Note that in the case $x=<a^1|b^1>$ and $y=<a^2|b^2>$ are  on the same unstable manifold, then there exists an $N>0$ such that  $a^1_j=a^2_j$, for any $j< -N$.  Therefore, when $\hat{A}: \hat{\Omega} \to \mathbb{R}$ is Holder and $(x,y)\in G$ then it is well defined
$$ \delta(x,y)= \Pi_{i=1}^\infty \frac{e^{\hat{A} (\hat{\sigma}^{-i} (x))}}{e^{\hat{A} (\hat{\sigma}^{-i} (y))}}= \Pi_{i=1}^\infty \frac{e^{\hat{A} (\hat{\sigma}^{-i} (<a^1|b^1>))}}{e^{\hat{A} (\hat{\sigma}^{-i} (<a^2|b^2>))}} . $$
Fix a certain $x_0=<a^0,b^0>$, then  the above can also be written as
$$ \delta(x,y)= \frac{e^{V(
x)}}{e^{V(y)}}= \frac{e^{V(<a^1|b^1>)}}{e^{V(<a^2|b^2>)}},$$
where
$$e^{V(<a|b>)}=  \Pi_{i=1}^\infty \frac{e^{\hat{A} (\hat{\sigma}^{-i} (<a|b>))}}{e^{\hat{A} (\hat{\sigma}^{-i} (<a^0|b^0>))}}     .  $$

Then, in this case $\delta$ is also of the form of example \ref{bib1}.

In this case, given any Haar system $\nu$ and any probability $m$,  Theorem \ref{gou} can be applied and we get examples of quasi-invariant probabilities.

\end{example}

\medskip

The next result has a strong similarity with the reasoning of  \cite{Man} and \cite{Seg}.

\medskip

\begin{proposition} \label{rod}  Given the modular function $\delta$ of example \ref{baker1}
consider the probability $M(d\,a,d\,b)=  W (b) \,d \,b\,d\,a\,  $ on $S^1 \times S^1$. Assume
$\nu^y$, $y=(a_0,b_0)$, is the Lebesgue probability $d b$ on the fiber $(a_0,b), 0\leq b <1$,
then, $M$ satisfies for all $f$

\begin{equation} \label{loip} \int \, \int f(s,y) \nu^y (ds) d M(y)= \int \, \int f(y,s) \delta^{-1}(y,s)\,\nu^y (ds) d M(y).
\end{equation}

\end{proposition}

{\bf Proof:}

We consider the equivalence relation:
given two points $z_1,z_2\in S^1 \times S^1$ they are related if the first coordinate is equal.

 In the case of example \ref{baker1} we take the a priori transverse function
$\nu^{z_1} (d\,b)= \nu^{a} (d\,b)$, $z_1=(a,\tilde b)$, constant equal to $d\,b$ in each fiber. This corresponds to the Lebesgue probability on the fiber.

For each pair $z_1=(a,b)$ and $z_2=(a,\,s)$, and $n \geq 0$,  the elements  $z_1^n, z_2^n$, $n \in \mathbb{N}$, such that $F^n (z_1^n)=z_1=(a,b)$ and  $F^n (z_2^n)=z_2=(a,s)$, are of the form
$z_1^n= (a^n, b^n),$ $ z_2^n=(a^n, s^n)$.

\bigskip

We define the cocycle
$$ \delta(z_1,z_2)=\Pi_{j=1}^\infty\,\, \,\frac{A ( z_1^n )}{A (z_2^n)}.$$

Fix a certain point $z_0=(a,c)$ and define $V$ by
$$ V(z_1)=\Pi_{j=1}^\infty\,\, \,\frac{A ( z_1^n )}{A (z_0^n )}.$$

Note that we can write
$$ \delta(z_1,z_2)= \frac{V(z_1)}{V(z_2)},$$
for such function $V$.

Remember that by notation $x_0$ is a point where $(0,x_0)$ and $(x_0,1)$ are intervals which are domains of injectivity of $T$.

\medskip

{\bf Remark:} Note the important point that if $x= (a,b)$ and $x'=(a',b)$, with $x_0\leq a \leq a'$, we get that
$b_n(x)=b_n(x').$ In the same way if $0\leq a \leq x_0$ we get that
$b_n(x)=b_n$. In this way the $b_n$ does not depends on $a$.

This means, there exists $W$ such that  we can write
$$\delta(z_1,z_2)=\delta^{-1}(\,(a,b),(a,s)\,)= Q(s,b)= \frac{W(s)}{W(b)},$$
where $b,s \in S^1$.

Condition (\ref{loip}) for $y$ of the form $y=(a,b)$ means for any $f$:
$$\int \, \int f(\,(a,b),(a,s)\,) \,\,\nu^a (d\,s) \,\,d M(a,b)= $$
$$\int \, \int f(\,(a,s),(a,b)\,) \,\,\delta^{-1}(\,(a,b),(a,s)\,)\,\,\nu^a (d\,s) \,\,d M(a,b)=
$$
$$\int \, \int f(\,(a,s),(a,b)\,) \,\,Q(s,b)\,\,\nu^a (d\,s) \,\,d M(a,b).
$$

Now, considering above  $f(\,(a,b),(a,s)\,) V(s)$ instead of $f(\,(a,b),(a,s)\,)$, we get the equivalent condition: for any $f$:

$$\int \, \int f(\,(a,b),(a,s)\,) \,W(s)\,d\,s \,\,d M(a,b)= $$
$$\int \, \int f(\,(a,s),(a,b)\,) \,\, W(s)\,\,d\,s \,\,d M(a,b).
$$

As $d M=
W (b) d \,b\,d a\,$ we get the alternative condition

$$\int \, \int f(\,(a,b),(a,s)\,) \,W(s)\,d\,s \,\,W (b) d \,b\,d\,a\,= $$
\begin{equation} \label{Rof} \int \, \int f(\,(a,s),(a,b)\,) \,\, W(s)\,\,d\,s \,\,W (b) d \,b\,d\,a\,,
\end{equation}
which is true because we can exchange the variables $b$ and $s$ on the first term above.

\qed

\begin{example} \label{gor} Consider the groupoid $G$ associated to the equivalence relation of example \ref{win1}.
In this case $x$ and $y$ are  on the same class when there exists an $N>0$ such that  $x_j=y_j$, for any $j \geq N$. Each class has a countable number of elements.

Consider a Holder potential $A: \overrightarrow{\Omega}\to \mathbb{R}$.

For $(x,y)\in G$ it is well defined
$$ \delta(x,y)= \Pi_{i=0}^\infty \frac{e^{A (\sigma^{i} (x))}}{e^{A (\sigma^{i} (y))}} . $$

 Consider the counting Haar system $\nu$ on each class.

We say  $f:G \to \mathbb{R}$ is admissible if for each class there exist  a finite number of non zero elements.

The quasi-invariant condition (\ref{kwe}) for the probability $M$ on $\overrightarrow{\Omega}$ means:
for any admissible integrable function  $f:G \to \mathbb{R}$  we have
\begin{equation} \label{kwe56}
 \sum_s\, \int f(s,x) d M(x)= \sum_s\, \int f(x,s) \,\Pi_{i=0}^\infty \frac{e^{A (\sigma^{i} (s))}}{e^{A (\sigma^{i} (x))}} \, d M(x).
\end{equation}

Suppose $B$ is such that $B= A + \log h - \log (g \circ \sigma) - c.$ This expression is called a coboundary equation for $A$ and $B$.  Under this assumption, as $x \sim s$, we get
$$\sum_s\, \int f(x,s) \,\Pi_{i=0}^\infty \frac{e^{B (\sigma^{i} (x))}}{e^{B (\sigma^{i} (s))}} \, d M(x)=$$
$$\sum_s\, \int f(x,s) \,\Pi_{i=0}^\infty \frac{e^{A (\sigma^{i} (x))}}{e^{A (\sigma^{i} (s))}} \frac{h(x)}{h(s)}\, d M(x).$$

Take $f(s,x)=g(s,x)\, h(x)$, then, as $M$ is quasi-invariant for $A$, we get that
\begin{equation} \label{kwe57}
 \sum_s\, \int g(x,s) \,\Pi_{i=0}^\infty \frac{e^{B (\sigma^{i} (x))}}{e^{B (\sigma^{i} (s))}} \, h(x)\, d M(x)= \sum_s\, \int g(s,x)\,h(x)\, d M(x).
\end{equation}

As $g(x,s)$ is a general function we get that $h(x) \, d M(x)$ is quasi-invariant for $B$.

Any Holder function $A$ is coboundary to a normalized Holder potential. In this way, if we characterize the quasi-invariant probability $M$ for  any given normalized potential $A$, then, we will be  able to determine, via the corresponding coboundary equation,
the quasi-invariant probability for any Holder potential.

\end{example}


\medskip




\begin{thebibliography}{99}

\bibitem{AHR}
Z. Afsar, A. Huef and I. Raeburn,  KMS states on $C^*$-algebras associated to local homeomorphisms, Internat. J. Math. 25 (2014), no. 8, 1450066, 28 pp

\bibitem{AnaR} C. Anatharaman-Delaroche, Ergodic Theory and Von Neumann algebras: an introduction, preprint Univ d'Orleans (France)

\bibitem{ADR} C. Anantharaman-Delaroche and J. Renault, Amenable groupoids, Monographs of L'Enseignement Math�matique, 36. L'Enseignement Math�matique, Geneva, 2000. 196 pp.

\bibitem{Araki} H. Araki,
On the
Equivalence of KMS and Gibbs Conditions for States of Quantum Lattice
Systems, Commun. Math. Phys, 35, 1--12 (1974)

\bibitem{Ba} T. Banakh, Direct Limit topologies and a characterization  of LF-spaces



\bibitem{BCLMS} A. Baraviera, L. M. Cioletti,  A. O. Lopes,
J. Mohr and R. R. Souza, On the general $XY$ Model: positive and zero temperature, selection and non-selection, \emph{Rev. in Math Phys}, Vol. 23, N. 10, pp 1063-1113 (2011).


\bibitem{Bis}  R. Bissacot and B. Kimura,
Gibbs Measures on Multidimensional Subshifts, preprint (2016)

\bibitem{Bla}
B. Blackadar, Operator algebras. Theory of $C^*$-algebras and von Neumann
algebras. Encyclopaedia of Mathematical Sciences, 122. Operator Algebras
and Non-commutative Geometry, III. Springer-Verlag, Berlin, 2006.

\bibitem{Bra} O. Bratteli and D. Robinson, Operator Algebras and Quantum Statistical Mechanics I and II, Springer Verlag.



  \bibitem{CDLS}
L. Cioletti, M. Denker, A. O. Lopes and M. Stadlbauer,
Spectral Properties of the Ruelle Operator for Product Type Potentials on Shift Spaces,
Journal of the London Mathematical Society, Volume 95, Issue 2, 684--704 (2017)


\bibitem{CL1}
L. Cioletti and A. O. Lopes.
Interactions, specifications, DLR probabilities and the Ruelle
operator in the one-dimensional lattice, Discrete and Cont. Dyn. Syst. - Series A, Vol 37, Number 12, 6139 -- 6152 (2017)


    \bibitem{Con} A. Connes, Sur la Theorie commutative de l'integration, Lect. Notes in Math. 725, Semminaire sur les Algebres d'Operateurs, Editor  P. de la Harpe,
19--143 (1979)

 \bibitem{Con1} A. Connes, Noncommutative Geometry, Academic Press (1994)


 \bibitem{Dea} V. Deaconu. Groupoids associated with endomorphisms. Trans. Amer. Math. Soc., 347(5):1779-1786,
 (1995).

    \bibitem{GiaII} G. DellAntonio,
Lectures on the Mathematics
of Quantum Mechanics II, Atlantis Press (2016)

 \bibitem{Dix} J. Diximier, Von Neumann Algebras, North Holland Publishing (1981)

  \bibitem{Exel1} R. Exel, A new look at the crossed-product of a $C^*$-algebra by an endomorphism. Ergod. Th. and Dynam. Sys. 23(6) (2003), 1733-1750.


    \bibitem{Exel2}  R. Exel, Crossed-products by finite index endomorphisms and KMS states. J. Funct. Anal. 199(1) (2003), 153-188.

     \bibitem{Exel3} R. Exel, KMS states for generalized gauge actions on Cuntz-Krieger algebras (An application of the Ruelle-Perron-Frobenius Theorem),  Bol. Soc. Brasil. Mat. (2004)

          \bibitem{Exel9}  R. Exel,     Inverse semigroups and combinatorial $C^*$-algebras,
Bull. Braz. Math. Soc. (N.S.), 39 (2008), 191--313.

\bibitem{EL1} R. Exel and A. Lopes,
$C^*$-algebras, approximately proper equivalence
relations and thermodynamic formalism, Ergod. Theo. and Dynam. Sys., 24, 1051-1082 (2004)

\bibitem{EL2} R. Exel and A. Lopes, $C^{*}$- Algebras and  Thermodynamic Formalism,
Sao Paulo Journal of Mathematical Sciences - Vol. 2, 1, 285--307 (2008)


 \bibitem{FM1} J. Feldman and C. Moore, Ergodic equivalence relations, cohomologies, von Neumann algebras I, TAMS 234 (1977) 289-359.


 \bibitem{FM2} J. Feldman and C. Moore, Ergodic equivalence relations, cohomologies, von Neumann algebras II, TAMS 234 (1977).

\bibitem{Hahn2} P. Hahn, The regular representations of measure groupoids. Trans. Amer. Math. Soc. 242 (1978), 35--72


\bibitem{HD}
N. T. A. Haydn and D. Ruelle,
Equivalence of Gibbs and Equilibrium states for homeomorphisms satisfying expansiveness and specification, Comm. in Math. Phys., 148, 155-167, 1992



\bibitem{Kas} D. Kastler, On Connes' Noncommutative Integration Theory, Comm. in Math. Physics, v 85, 99--120 (1982)


\bibitem{ManS}
M. Kessebohmer, M. Stadlbauer and B. Stratmann, Lyapunov spectra for KMS states on Cuntz-Krieger algebras. Math. Z. 256 (2007), no. 4, 871--893.



 \bibitem{KumRen} A. Kumjian and J. Renault, KMS states on   $ C^*$-Algebras associated to expansive maps, Proc. AMS Vol. 134, No. 7, 2067--2078 (2006)


\bibitem{Kum1}
A. Kishimoto and A.  Kumjian,  Simple stably projectionless $C^*$-algebras arising as crossed products. Canad. J. Math. 48 (1996), no. 5, 980--996.


\bibitem{Lan} N. Landsman, Lecture Notes on $C^*$-Algebras and Quantum Mechanics, Univ. of Amsterdam (1998)

\bibitem{LY} F. Ledrappier and L.-S. Young,
The Metric Entropy of Diffeomorphisms: Part I: Characterization of measures satisfying Pesin's Entropy Formula,
Annals of Mathematics,
Vol. 122, No. 3,  509--539 (1985)

\bibitem{LM2} A. O.  Lopes  and G. Mantovani,
The KMS Condition for the homoclinic equivalence relation and Gibbs probabilities, to appear in Sao Paulo Jour. of Math. Scien.

\bibitem{LO} A. O.  Lopes  and E. Oliveira, Continuous groupoids on the symbolic space, quasi-invariant probabilities for Haar systems and the Haar-Ruelle operator, to appear in Bull of the Braz. Math. Soc.




\bibitem{LMe} A. O. Lopes  and J. K. Mengue, Thermodynamic Formalism for Haar systems in Noncommutative Integration: transverse functions and entropy of  transverse measures, in preparation





\bibitem{LMMS2} A. Lopes, J. K. Mengue, J. Mohr and R. R. Souza,
Entropy and variational principle for one-dimensional Lattice systems with a general a-priori probability: positive and zero temperature, Erg. Theory and Dyn Systems, 35 (6), 1925--1961 (2015)

\bibitem{Meye}
T. Meyerovitch, Gibbs and equilibrium measures for some families of subshifts,
Ergod. Th. and Dynam. Sys., 33, 934�-953, (2013)

\bibitem{Man} G. Mantovani, Teoria n\~ao comutativa de integra\c c\~ao e din\^amica hiperb\'olica, Disserta\c c\~ao de Mestrado,  ICMC-USP - S\~ao Carlos - Brasil (2013)



\bibitem{Ou} M. O'Uchi, Measurable groupoids and associated von Neumann algebras, Ehime University Notes (1984)


\bibitem{Pana}  J. Panangaden,
Energy Full Counting Statistics in Return-to-Equilibrium
Jane Panangaden, thesis McGill Univ (2016)


\bibitem{PP}
W. Parry and M.  Pollicott. Zeta functions and the periodic
orbit structure of hyperbolic dynamics. \emph{Ast\'erisque}
Vol {187-188} (1990).


\bibitem{PY}
M.  Pollicott and M. Yuri, Dynamical systems and Ergodic Theory, Cambrige Press, 1998


\bibitem{Ped}
G. K. Pedersen, $C^*$-Algebras and and their automorphism groups, Acad Press, (1979)

\bibitem{Put} I. Putnam,
Lecture Notes on Smale Spaces, Univ. of Victoria - Canada (2015)

\bibitem{PutJ} I. Putnam and J. Spielberg, The Structure of $C^*$-Algebras associated with hyperbolic dynamical systems, J. Funct. Anal. 163 (1999), no. 2, 279--299

\bibitem{Ren0}    J. Renault, A Groupoid approach to $C^*$-algebras, Lecture Notes in Mathematics 793, Springer-Verlag, (1980)



\bibitem{Ren1}
J. Renault,
AF equivalence relations and their cocycles, Operator algebras and mathematical physics (Constanca, 2001), 365--377, Theta, Bucharest, (2003).

\bibitem{Ren1.5}
J. Renault,
The Radon-Nikodym problem for appoximately proper equivalence relations.
Ergodic Theory Dynam. Systems 25 (2005), no. 5, 1643--1672.

\bibitem{Ren2}
J. Renault, $C^*$-Algebras and Dynamical Systems, XXVII Coloquio Bras. de Matematica - IMPA (2009)


\bibitem{Ru1} D. Ruelle, Noncommutative algebras for hyperbolic diffeomorphisms, Inv. Math. 93, 1-13 (1988)



\bibitem{Seg} J. Segert, Hyperbolic Dynamical Systems and the Noncommutative Theory of Connes, PhD Thesis, Department of Physics, Princeton University (1987)


\bibitem{Bak}
L. Smith,
Chaos: A Very Short Introduction, Cambridge Univ. Press

 \bibitem{Tak} M. Takesaki,   Tomita's Theory
of Modular Hilbert Algebras
and its Applications,  Lecture Notes in Mathematics 128, Springer-Verlag, (1970)


\bibitem{Thom1}
K. Thomsen,  The homoclinic and heteroclinic $C^*$-algebras of a generalized one-dimensional solenoid, Ergodic Theory Dynam. Systems 30, no. 1, 263--308 (2010).


\bibitem{Thom2}
K. Thomsen, $C^*$-algebras of homoclinic and heteroclinic structure in expansive dynamics. Mem. Amer. Math. Soc. 206, no. 970 (2010)

\bibitem{Thom3}
K. Thomsen, On the $C^*$-algebra of a locally injective surjection and its KMS states. Comm. Math. Phys. 302, no. 2, 403--423 (2011)

\bibitem{Thom4} K. Thomsen, Phase transitions on $O_2$, Comm. Math. Phys, 349 (2), 481-492 (2017)

\bibitem{Thom5}
K. Thomsen, KMS-states and conformal measures, Comm. Math. Phys. 316 (2012),
615-640

\bibitem{Thom6} K. Thomsen, KMS weights on groupoid and graph $C^*$-algebras, J. Func. Analysis 266
(2014), 2959-2998.

\bibitem{Thom7}
K. Thomsen, KMS weights on graph $C^*$-algebras II. Factor types and ground states, Arxiv



\end{thebibliography}
\end{document}